\documentclass[11pt, reqno]{amsart}
\usepackage[margin=1in]{geometry}

\usepackage{amsmath,amssymb,amsthm,mathtools}
\usepackage{enumitem}
\usepackage{charter}
\usepackage{newtxmath}

\allowdisplaybreaks

\usepackage{xcolor} %
\definecolor{darkred}{rgb}{0.5,0,0} %
\definecolor{darkblue}{rgb}{0,0,0.5} %

\usepackage[colorlinks,linkcolor=black,citecolor=darkblue, pagebackref=false,urlcolor=darkblue,hypertexnames=true]{hyperref}

\usepackage{euflag}

\usepackage{microtype}

\usepackage{needspace}
\usepackage{longtable,booktabs,array}

\makeatletter
\renewcommand{\l@section}{\@tocline{1}{6pt}{1pc}{}{}}
\renewcommand{\l@subsection}{\@tocline{2}{3pt}{1pc}{5pc}{}}
\makeatother

\newtheorem{theorem}{Theorem}[section]

\newtheorem{corollary}[theorem]{Corollary}
\newtheorem{lemma}[theorem]{Lemma}

\newtheorem{proposition}[theorem]{Proposition}

\newtheorem*{theorem*}{Theorem}

\theoremstyle{definition}

\theoremstyle{definition}
\newtheorem{remark}{Remark}[section]

\theoremstyle{plain}
\newtheorem{introthm}{Theorem}

\newcommand{\vd}{{\tt{d}}}
\newcommand{\ve}{{\tt{e}}}
\newcommand{\vq}{{\tt{q}}}

\numberwithin{equation}{section}

\newcommand{\R}{\mathbb R}
\newcommand{\C}{\mathbb C}
\newcommand{\Hh}{\mathcal H}
\newcommand{\B}{\mathcal B}
\newcommand{\SOS}{\operatorname{SOS}}
\newcommand{\diag}{\operatorname{diag}}
\newcommand{\adj}{\operatorname{adj}}
\newcommand{\col}{\operatorname{col}}
\newcommand{\M}{\mathcal M}

\makeatletter
\renewcommand{\subsubsection}{\@startsection{subsubsection}{3}%
  \z@{1.5ex plus .5ex minus .2ex}%
  {0.5ex plus .2ex}{\normalfont\normalsize\bfseries}}
\makeatother

\title[Operator Fejer--Riesz Theorems]{Operator Fej\'er--Riesz Theorems with Gaps and Degree Bounds}

\author[A.\ Jindal]{Abhay Jindal}
\address{Abhay Jindal, Faculty of Mathematics and Physics, University of Ljubljana, Slovenia}
\email{abhay.jindal@fmf.uni-lj.si}
 \thanks{AJ was supported by the Slovenian Research Agency program P1-0222 and
grants J1-50002, N1-0536. This work was performed
within the project COMPUTE, funded within the QuantERA II Programme, which has
received funding from the European Union's Horizon 2020 research and innovation
programme under Grant Agreement No.~101017733 {\normalsize\protect\euflag}.}

\author[A.\ Zalar]{Aljaž Zalar}
\address{Aljaž Zalar, University of Ljubljana, Faculty of Computer and Information Science \&
Faculty of Mathematics and Physics, \& Institute of Mathematics, Physics and Mechanics,
Ljubljana, Slovenia.}
\email{aljaz.zalar@fri.uni-lj.si}
\thanks{AZ was supported by the Slovenian Research Agency  research core funding No. P1-0288 and grants J1-50002, J1-60011, J1-70017.}

\date{}
\begin{document}
\begin{abstract}
Let $K$ be a nonempty compact semialgebraic subset of the real line, and let $S$ be a saturated scalar description of $K$. We prove that every hermitian polynomial $F\in\mathcal B(\mathcal H)[x]$ positive semidefinite on $K$ admits a certificate in the operator quadratic module generated by $S$, with weighted degree at most $\deg F+c_S$. Here, $c_S$ depends only on $S$ and is independent of both $F$ and the coefficient Hilbert space $\mathcal H$, which may be infinite-dimensional.

\smallskip

In contrast to the scalar case, natural matrix preorderings admit no uniform degree bound over compact semialgebraic subsets of the real line. We construct positive semidefinite quadratic $2\times2$ matrix polynomials on unions of two intervals whose minimum certificate degrees tend to infinity as one interval shrinks to a point, while the gap remains bounded below. Conversely, uniform degree bounds hold for interval unions in a fixed compact interval when the component and gap lengths are bounded below by positive constants.

\smallskip

Our proofs develop two techniques. The first is a quadratic reduction principle that reduces the problem of obtaining degree bounds for arbitrary-degree operator polynomials to the quadratic case, with a fixed additive increase in degree. The second uses a rank-one symmetric matrix pencil to transfer positivity on unions of intervals to $[0,1]$, where the operator Fej\'er--Riesz theorem provides certificates. We also obtain operator-valued Fej\'er--Riesz factorizations on unions of circular arcs, rational certificates with a fixed denominator on unbounded semialgebraic sets, and a characterization of the sets whose natural operator preorderings admit degree-preserving certificates.
\end{abstract}

\keywords{Positivstellensatz, operator polynomials, matrix polynomials, sums of squares, quadratic module, preordering, degree bounds, Fejér–Riesz factorization, semialgebraic sets, real algebraic geometry}
\subjclass[2020]{Primary 47A56, 13J30;
Secondary 14P10, 47A68, 90C23}

\maketitle

\clearpage

\tableofcontents

\clearpage

\section{Introduction}\label{sec:intro}
\subsection{Background and motivation}
Positivity and sums of squares play a central role in real
algebraic geometry and operator theory. The search for
algebraic certificates of polynomial nonnegativity goes
back to Hilbert's $17^{\rm th}$ problem and is closely
related to moment theory and the Positivstellens\"atze
of Schm\"udgen and Putinar \cite{Schmudgen1991,Putinar1993};
for classical results and modern treatment see \cite{BCR98,Mar08,Sc24}.
Two fundamental questions arise in this context:
when does a nonnegative polynomial admit a sums-of-squares
certificate, and how can one bound the degree of such
a certificate?

\medskip

The operator Fej\'er--Riesz theorem~\cite{Rosenblum1968,DR} provides a sharp answer to both questions for trigonometric polynomials in one variable. If a trigonometric polynomial
$P(e^{i\theta})=\sum_{k=-\vd}^{\vd}A_k e^{i k\theta}$
with coefficients $A_k\in\B(\Hh)$ is positive
semidefinite for every $\theta$, then
\[
P=Q^*Q
\]
for an analytic polynomial
$Q(e^{i\theta})=\sum_{k=0}^{\vd}B_k e^{i k\theta}$
with $B_k\in\B(\Hh)$. Fej\'er--Riesz factorization plays an important role in the operator-algebraic theory of Toeplitz operator systems \cite{CvS21, Farenick2021, DuhanJindal}. Related developments include multivariable factorization \cite{GW05,Dritschel,DW05, Dritschel2025, DKMV26} and noncommutative positivity and factorization results \cite{He,Mc}.

\medskip

While the Fej\'er--Riesz theorem concerns positivity on
the entire unit circle (or, in its real-line formulation,
on the entire real line), applications often require
positivity on prescribed subsets of $\R$. Such problems arise in systems and control theory, where finite-frequency generalizations of the Kalman--Yakubovich--Popov lemma address positivity on restricted frequency intervals \cite{IH05,PIH14}.

\medskip

We study hermitian polynomials with bounded-operator coefficients that are positive semidefinite on arbitrary compact semialgebraic subsets of the real line, that is, unions of intervals and isolated points. We ask whether such polynomials admit exact positivity certificates in the corresponding operator quadratic modules, whether the certificate degrees can be bounded independently of the coefficient Hilbert space, and how these bounds depend on the underlying semialgebraic set.  

\subsection{Notation}\label{ssec:Notation}
Let $\Hh$ be a real or complex Hilbert space, and let
$\B(\Hh)[x]$ denote the polynomials in a real indeterminate
$x$ with coefficients in the bounded operators on $\Hh$.
For a nonzero polynomial $F(x)=\sum_{k=0}^{\vd} A_kx^k$,
with $A_{\vd}\ne0$, we define $\deg F= \vd$ and set
$\deg0=-\infty$. We write
\[
F^*(x)=\sum_{k=0}^{\vd} A_k^*x^k
\]
and call $F$ \emph{hermitian} if $F^*=F$, or equivalently,
if $A_k^*=A_k$ for every $k$. Finally, we write
$F\succeq0$ on $K$ if $F(x)$ is positive semidefinite
for every $x\in K$.

\medskip

For a complex Hilbert space $\Hh$, let $\SOS_{\Hh}=\{G^*G:G\in\B(\Hh)[x]\}$. For a real Hilbert space, we use $\SOS_{\Hh}=\{G^*G+H^*H:G,H\in\B(\Hh)[x]\}$. The fact that these definitions agree with the usual finite sums of hermitian squares follows from the factorization results recalled in Section~\ref{sec:tools} and the real descent argument in Appendix~\ref{app:real}. We also set $\SOS_{\C}=\{f^*f:f\in\C[x]\} \subset \R[x]$.

\medskip

Let $S=\{g_1,\ldots,g_s\}\subseteq\R[x]$ be finite, and put $g_0=1$. Define
\begin{align*}
K_S&=\{x\in\R \ \colon \ g_j(x)\ge0,\ 1\le j\le s\},\\
M(S)&=\biggl\{\sum_{j=0}^s g_j\sigma_j \ \colon\ \sigma_j\in\SOS_{\C}\biggr\},\\
\M_{\Hh}(S)&=\biggl\{\sum_{j=0}^s g_j\Sigma_j \ \colon\ \Sigma_j\in\SOS_{\Hh}\biggr\}.
\end{align*}
We call $M(S)$ the \emph{quadratic module generated by $S$}
and $\mathcal M_{\mathcal H}(S)$ the
\emph{operator quadratic module generated by $S$}. The module $M(S)$ is \emph{saturated}
if it contains every real polynomial nonnegative on $K_S$; in this case, $S$ is called a \emph{saturated scalar description} of $K_S$.

\medskip

The \emph{weighted degree} of an expression $\sum_j g_j\Sigma_j$ is the maximum of the degrees of its summands. We write $\M_{\Hh}(S)_{\le n}$ for the polynomials admitting such a certificate of weighted degree at most $n$. The \emph{operator preordering} generated by $S$ is
\[
\mathcal T_{\Hh}(S)=\biggl\{\sum_{\alpha\in\{0,1\}^s}g^\alpha\Sigma_\alpha \ \colon\  \Sigma_\alpha\in\SOS_{\Hh}\biggr\},\qquad
 g^\alpha=\prod_{j=1}^s g_j^{\alpha_j}.
\]
The corresponding truncated preordering, with weighted degree at most $n$, is denoted by $\mathcal T_{\Hh}(S)_{\le n}$.

\medskip

These definitions also apply when $\Hh=\C^m$, in which case we speak of matrix quadratic modules and matrix preorderings. For $n\ge0$, we write $\R[x]_{\le n}$ for the polynomials of degree at most $n$, and use $\mathbb N=\{1,2,\ldots\}$.

\subsection{Main results}

We are now ready to state our first result, which is the main contribution of this paper.

\begin{introthm}\label{thmA}\label{thm:main}
Let $S=\{g_1,\ldots,g_s\}\subseteq\R[x]$ be finite,
with $K_S$ nonempty and compact.
Suppose that the scalar quadratic module $M(S)$ is saturated.
Then, for every real or complex Hilbert space $\Hh$,
every nonnegative integer $\vd$, and every hermitian
$F\in\B(\Hh)[x]$ with $\deg F\le\vd$,
\[
F\succeq0\text{ on }K_S
\quad\Longleftrightarrow\quad
F=\sum_{j=0}^{s}g_j
\sum_{i=1}^{r}F_{j,i}^*F_{j,i},
\]
where $g_0=1$, $F_{j,i}\in\B(\Hh)[x]$, and
$r=1$ if $\Hh$ is complex and $r=2$ if $\Hh$ is real.
Moreover, there exists a nonnegative integer $c_S$,
depending only on $S$ (and independent of $\Hh, \vd, F$), such that the factors can be chosen
to satisfy
\[
\deg(g_jF_{j,i}^*F_{j,i})
\le \vd+c_S
\qquad (0\le j\le s,\ 1\le i\le r).
\]
\end{introthm}

\smallskip

\begin{remark}
Several remarks related to Theorem~\ref{thmA} are in order. 
\begin{enumerate}[leftmargin=*,font=\normalfont,itemsep=3pt]
\item[(a)] Every nonempty compact semialgebraic set $K\subseteq\R$ admits a saturated scalar description, and Scheiderer characterized all such descriptions \cite{Scheiderer2003,Scheiderer2005Distinguished}; see Theorem~\ref{thm:scalar}. For natural descriptions, the scalar preorderings are saturated \cite{KM} and admit certificates of weighted degree at most $\deg f$ \cite{KMS}.

\item[(b)] For a saturated scalar description $S$ of a compact semialgebraic set, the second named author \cite{Zalar} proved a non-quantitative version of Theorem~\ref{thmA} for matrix polynomials. Sun and the second named author \cite{SunZalar} obtained degree bounds, independent of matrix size, when the set consists of a compact interval and an isolated point.

\item[(c)] The degree bound in Theorem~\ref{thmA}
is optimal up to a fixed additive constant.
Indeed, fix $a<\min K_S$.
For every $\vd\ge0$, the polynomial
$(x-a)^{\vd}$ is nonnegative on $K_S$,
and any certificate for it must have
weighted degree at least $\vd$.
Thus every uniform degree-bound function
must be bounded below by $\vd$, whereas
Theorem~\ref{thmA} provides the upper
bound $\vd+c_S$.

\item[(d)] If the generators have real algebraic coefficients, the proof also gives a way to compute an upper bound on the constant $c_{S}$ by exact semidefinite feasibility tests; see Subsection~\ref{subsec:degree-constants}.
\item[(e)] For compact semialgebraic sets of dimension at least two,
the weighted degree of positivity certificates generally
cannot be bounded solely in terms of the polynomial degree
and the fixed defining generators, even in the scalar
case \cite{Scheiderer2005Nonexist}.

\item[(f)] The classical Positivstellens\"atze of Schm\"udgen and Putinar provide certificates for strictly positive polynomials under compactness or Archimedean assumptions \cite{Schmudgen1991,Putinar1993}. The degree bounds typically depend on the polynomial degree, its coefficients, and its minimum on the underlying set \cite{Schweighofer2004,NieSchweighofer2007,
BaldiMourrain2023,BaldiSlot2024}.
Related matrix and operator Positivstellensätze in several variables have been established under additional hypotheses \cite{AmbrozieVasilescu2003,SchererHol2006,KlepSchweighofer2010,Cimpric2011Operator,Cimpric2011Strict}. Other developments include Positivstellensätze for matrix algebras \cite{SavchukSchmudgen2012,Cimpric2012Matrix} and quantitative degree bounds for matrix polynomials \cite{Huang2025}.
\end{enumerate}
\end{remark}

Although Theorem~\ref{thmA} gives degree bounds independent of the coefficient Hilbert space, the additive constant depends on the scalar description and thus on the underlying semialgebraic set. Our second result shows that this dependence cannot be avoided even for natural matrix preorderings and quadratic $2\times2$ polynomials.

\medskip

Fix $\ell \in \mathbb N$. For a finite union of nondegenerate compact intervals
\[
K=\bigcup_{j=1}^{\ell}[a_j,b_j],
\qquad a_1<b_1<\cdots<a_\ell<b_\ell,
\]
we write
\[
S_K=\{x-a_1,b_\ell-x\}
\cup\{(x-b_j)(x-a_{j+1}) \ :\ 1\le j<\ell\}
\]
for its natural scalar description.

 \begin{introthm}\label{thmB}
    For each $0<\varepsilon<1/2$, there exist sets
    $K_{\varepsilon}\subseteq[0,1]$ and hermitian matrix
    polynomials $F_{\varepsilon}\in\C[x]^{2\times2}$
    of degree at most $2$ such that the following hold:
 \begin{enumerate}[label=(\roman*),
                  font=\normalfont,
                  itemsep=3pt]
        \item Each $K_{\varepsilon}$ is a union of two disjoint
        compact intervals, one of which shrinks to a point
        as $\varepsilon\to0$, while the gap between intervals
        has length bounded below independently of $\varepsilon$.
    
        \item $F_{\varepsilon}\succeq0$ on $K_{\varepsilon}$.
    
        \item If $S_{K_\varepsilon}$ is the natural description of
        $K_{\varepsilon}$, then
        \[
        \min\bigl\{n \ \colon\ F_{\varepsilon}\in
        \mathcal T_{\C^2}(S_{K_\varepsilon})_{\le n}\bigr\}
        \longrightarrow\infty
        \qquad(\varepsilon\to0).
        \]
    \end{enumerate}
    \end{introthm}

This contrasts with the scalar case, where natural scalar preorderings admit certificates of weighted degree at most $\deg f$, independently of the underlying set \cite{KMS}. 

\medskip

Our third result shows that uniform additive degree bounds can nevertheless be established for operator-valued polynomials when interval lengths and gap lengths are bounded away from zero.

\medskip

For the interval union $K$ above, set
\[
\Delta_K(x)=-\prod_{j=1}^{\ell}(x-a_j)(x-b_j).
\]

\begin{introthm}\label{thmC}
Fix a compact interval $J\subset\R$ and $\mu,\omega>0$. Let $\mathcal K(J,\mu,\omega)$ be the family of all sets
\[
K=\bigcup_{j=1}^{\ell}[a_j,b_j]\subseteq J,\qquad \ell\ge1,
\]
where $a_1<b_1<\cdots<a_\ell<b_\ell$, each interval has length at least $\mu$, and each gap between consecutive intervals has length at least $\omega$. There exists a nonnegative integer $C$, depending only on $J,\mu,\omega$, such that for every $K\in\mathcal K(J,\mu,\omega)$, every real or complex Hilbert space $\Hh$, every nonnegative integer $\vd$, and every hermitian $F\in\B(\Hh)[x]$ with $\deg F\le\vd$,
\[
F\succeq0\text{ on }K
\quad\Longleftrightarrow\quad
F\in\M_{\Hh}(\{\Delta_K\})_{\le\vd+C}
\quad\Longleftrightarrow\quad
F\in\mathcal T_{\Hh}(S_K)_{\le\vd+C}.
\]
\end{introthm}

\smallskip

Theorem~\ref{thmB} shows that the lower bound on interval lengths in Theorem~\ref{thmC} is essential, even when the gap lengths remain bounded below.

\subsection{Further consequences}

Theorem~\ref{thmA} and its intermediate results have several consequences beyond
compact semialgebraic subsets of the real line.

\medskip

Let $\mathbb T=\{z\in\C:|z|=1\}$ denote the unit circle.
For an operator Laurent polynomial
$F(z)=\sum_{k=-m}^n A_kz^k$, we put
\[
F^*(z)=\sum_{k=-m}^n A_k^*z^{-k},
\qquad
\deg_L F=\max\{|k|:A_k\ne0\},
\]
and set $\deg_L0=-\infty$.
On $\mathbb T$, the value $F^*(z)$ is the adjoint of $F(z)$.
We call $F$ \emph{hermitian} if $F^*=F$, equivalently,
$A_{-k}=A_k^*$ for every $k$.

\medskip

 The notion of a saturated scalar description of a finite union of closed arcs of $\mathbb T$ is defined in Subsection~\ref{subsec:circle}.

\medskip

We obtain an operator-valued Fej\'er--Riesz
theorem on finite unions of closed arcs of the unit
circle, including isolated points. It
 extends \cite[Theorem~2.5.4]{ZalarThesis} to operator
coefficients and provides uniform degree bounds. 

\begin{corollary}
\label{cor:circle}
Let $E\subseteq\mathbb T$ be a nonempty finite union of arcs,
and let
$\mathcal S=\{g_1,\ldots,g_s\}\subseteq\C[z,z^{-1}]$
be a saturated description of $E$.
Put $g_0=1$.
Then, for every complex Hilbert space $\Hh$,
every nonnegative integer $\vd$, and every hermitian
$F\in\B(\Hh)[z,z^{-1}]$ with $\deg_L F\le \vd$,
\[
F\succeq0\text{ on }E
\quad\Longleftrightarrow\quad
F=\sum_{j=0}^{s}g_jQ_j^*Q_j,
\]
where $Q_j\in\B(\Hh)[z]$.
Moreover, there exists a nonnegative integer
$c_{\mathcal S}$, depending only on $\mathcal S$ (and independent of $\Hh, \vd, F$),
such that the factors can be chosen to satisfy
\[
\deg_L(g_jQ_j^*Q_j)\le \vd+c_{\mathcal S}
\qquad (0\le j\le s).
\]
\end{corollary}

Let $K\subseteq\R$ be a nonempty closed semialgebraic
set, written as
\[
K=\bigcup_{j=1}^{\ell}([a_j,b_j]\cap\R),
\]
where $a_j,b_j\in\R\cup\{-\infty,+\infty\}$ and
$ a_1\le b_1<a_2\le b_2<\cdots<a_\ell\le b_\ell. $
The \emph{natural description} of $K$ is
\[
S_K=E_1\cup E_2\cup E_3,
\]
where
\[
\begin{aligned}
E_1&=
\begin{cases}
\{x-a_1\},&a_1>-\infty,\\
\varnothing,&a_1=-\infty,
\end{cases}\\[2mm]
E_2&=\{(x-b_j)(x-a_{j+1}):1\le j<\ell\},\\[2mm]
E_3&=
\begin{cases}
\{b_\ell-x\},&b_\ell<+\infty,\\
\varnothing,&b_\ell=+\infty.
\end{cases}
\end{aligned}
\]

For unbounded semialgebraic sets, certificates in quadratic modules automatically preserve degree, since their highest-degree terms cannot cancel. Thus, Corollary~\ref{cor:degree-preservation} also characterizes the unbounded sets for which every positive semidefinite operator polynomial admits a certificate in the natural preordering. For other unbounded sets, such certificates are not available in general. The following result shows that allowing a fixed denominator overcomes this obstruction.

\begin{corollary}
\label{cor:unbounded}
Let $K\subseteq\R$ be an unbounded, closed semialgebraic set, and let
$S=\{g_1,\ldots,g_s\}\subseteq\R[x]$ be its natural
description. Put $g_0=1$.
There exists a nonnegative integer $c_S$, depending only
on $S$, such that, for every real or complex Hilbert
space $\Hh$, every nonnegative integer $\vd$, and every
hermitian $F\in\B(\Hh)[x]$ with $\deg F\le \vd$,
\[
F\succeq0\text{ on }K
\quad\Longleftrightarrow\quad
(1+x^2)^{c_S}F
 =\sum_{j=0}^{s}g_j\sum_{i=1}^{r}F_{j,i}^*F_{j,i},
\]
where $F_{j,i}\in\B(\Hh)[x]$, and $r=1$ if $\Hh$
is complex and $r=2$ if $\Hh$ is real.
Moreover, the factors can be chosen to satisfy
\[
\deg(g_j F_{j,i}^*F_{j,i})
 \le \vd + 2c_S
\qquad
(0\le j\le s,\ 1\le i \le r).
\]
If $K=\R$, one may take $S=\varnothing$ and $c_S=0$.
\end{corollary}

Finally, we characterize the closed semialgebraic
sets for which natural operator preorderings
admit certificates without any degree increment.

\begin{corollary}\label{cor:degree-preservation}
Let $K\subseteq\R$ be nonempty, closed and semialgebraic,
and let $S_K$ be its natural description.
Fix a real or complex Hilbert space $\Hh$ with
$\dim\Hh\ge2$.
Every nonzero hermitian $F\in\B(\Hh)[x]$ positive
semidefinite on $K$ belongs to
\[
\mathcal T_{\Hh}(S_K)_{\le\deg F}
\]
if and only if $K$ is one of the following:
\begin{enumerate}[label=(\roman*),
                  font=\normalfont,
                  itemsep=3pt]
\item a set of at most three points;
\item a nondegenerate compact interval together with
at most one additional isolated point;
\item $\R$, or a union of one or two disjoint unbounded intervals
together with at most one additional isolated point.
\end{enumerate}
\end{corollary}

\subsection{Proof overview and organization}
We give an overview of the proof of Theorem~\ref{thmA}. The starting point of our proof is the quadratic reduction
principle (Theorem~\ref{thm:quadratic-reduction}).
It shows that any degree bound for quadratic operator
polynomials that is independent of the coefficient Hilbert
space extends to a bound of the form $\vd+c_S$ for
polynomials of degree $\vd$.
The idea is to
represent an arbitrary-degree operator polynomial as a
compression of a quadratic polynomial on an enlarged
Hilbert space. A suitable perturbation makes the lifted
polynomial positive semidefinite, and the closedness of
truncated operator quadratic modules allows us to pass
to the limit.

\medskip

It remains to establish degree bounds for quadratic operator polynomials that are independent of the coefficient
Hilbert space. For
unions of nondegenerate compact intervals, we construct
a symmetric matrix pencil by a rank-one perturbation of
a diagonal matrix. As the parameter varies over $[0,1]$,
the eigenvalues of this pencil cover the prescribed
intervals. This allows us to transfer positivity to
$[0,1]$ and apply the operator Fej\'er--Riesz theorem.
However, pulling the resulting certificate back to the
original variable introduces a scalar multiplier that
vanishes in the gaps between the intervals. To overcome
this difficulty, we characterize the polynomial congruences
that preserve the quadratic degree and use them to obtain
a strictly positive multiplier. Scalar saturation then
allows us to remove this multiplier.
Both the multiplier and the scalar identity used to
remove it can be chosen to depend only on the underlying interval union and not on the operator polynomial.
This differs from the construction in \cite{Zalar},
where the auxiliary multiplier depends on the
particular matrix polynomial. As a consequence,
the degree increase in our construction can be
bounded solely in terms of the underlying interval
union.

\medskip

We handle isolated points using Douglas' factorization lemma,
polynomial interpolation, and certificates for auxiliary matrix
polynomials. Both the auxiliary polynomials and their certificates
depend only on the underlying set.
Combining
these results with the quadratic reduction principle
completes the proof of Theorem~\ref{thmA}.

\medskip

In Section~\ref{sec:tools} we collect the preliminary results.
Section~\ref{sec:reduction} establishes the quadratic
reduction principle and treats finite sets.
In Section~\ref{sec:intervals} we construct quadratic
certificates on unions of nondegenerate intervals,
while Section~\ref{sec:points} handles isolated
points and completes the proof of Theorem~\ref{thmA}.
In Sections~\ref{sec:nonuniform} and~\ref{sec:uniform}
we prove Theorems~\ref{thmB} and~\ref{thmC}, respectively.
Finally, Section~\ref{sec:further} proves the
corollaries stated in the introduction.

\section{Preliminaries}\label{sec:tools}
By Lemma~\ref{lem:real-descent}, it suffices to establish Theorems~\ref{thmA} and~\ref{thmC} and Corollaries~\ref{cor:unbounded} and~\ref{cor:degree-preservation} for complex Hilbert spaces.
We therefore work with complex Hilbert spaces throughout Sections~\ref{sec:tools}--\ref{sec:further}; the real case follows with the same weighted degree bounds.

\medskip

The following two consequences of the operator
Fej\'er--Riesz theorem \cite{DR,Rosenblum1968}
are well known. Proofs for real Hilbert spaces
can be found in \cite{CZ}; the arguments for
complex Hilbert spaces are essentially the same.

\begin{proposition}\label{prop:factorization}
Let $\vd$ be a nonnegative integer, and let
$F\in\B(\Hh)[x]$ be hermitian with
$\deg F\le2\vd$.
\begin{enumerate}[label=(\roman*),font=\normalfont, itemsep=3pt]
\item\label{item:global-factorization}
If $F\succeq0$ on $\R$, then
$F=G^*G$ for some $G\in\B(\Hh)[x]$
with $\deg G\le\vd$.

\item\label{item:interval-factorization}
If $F\succeq0$ on $[0,1]$, then there exist
$G,H\in\B(\Hh)[x]$ such that
\begin{equation}\label{eq:interval-factor}
F=G^*G+x(1-x)H^*H,
\qquad
\deg G\le\vd,\quad \deg H\le\vd-1,
\end{equation}
where $H=0$ if $\vd=0$.
\end{enumerate}
\end{proposition}

Fix $\ell\in\mathbb N$. Recall that every nonempty compact semialgebraic
subset $K\subseteq\R$ can be written as
\begin{equation}\label{def:semialgebraic}
K=\bigcup_{j=1}^{\ell}[a_j,b_j],
\qquad
a_j\le b_j,\quad b_{j-1}<a_j
\quad (2\le j\le \ell),
\end{equation}
where $a_1\le b_1$. In particular, an interval
$[a_j,b_j]$ with $a_j=b_j$ represents an isolated
point. 
The following characterization of saturation
follows from \cite[Theorem~5.17]{Scheiderer2003} and
\cite[Corollary 4.4]{Scheiderer2005Distinguished}.

\begin{theorem}
\label{thm:scalar}
Let $S=\{g_1,\ldots,g_s\}\subseteq\R[x]$ be finite
such that $K_S$ is nonempty and compact, written as in
\eqref{def:semialgebraic}.
Then $M(S)$ is saturated if and only if
the following conditions hold:
\begin{enumerate}[label=\textup{(\roman*)},itemsep=3pt]
\item For each $j=1,\ldots,\ell$, there exists
$i\in\{1,\ldots,s\}$ such that
    $g_i(a_j)=0$,
    $g_i'(a_j)>0.$
\item For each $j=1,\ldots,\ell$, there exists
$i\in\{1,\ldots,s\}$ such that
    $g_i(b_j)=0$,
    $g_i'(b_j)<0.$
\end{enumerate}
\end{theorem}

Scalar saturation on a compact subset of $\R$ extends to matrices
of every finite size by \cite[Theorem~2.1]{Zalar}. Theorem~\ref{thm:scalar} yields the scalar identity \eqref{eq:scalar},
which will allow us to cancel powers of a fixed strictly positive
polynomial in operator positivity certificates.

\begin{lemma}\label{lem:scalar-identity}
Let $\Delta\in\R[x]$ be such that
$K:=\{x\in\R:\Delta(x)\ge0\}$ is nonempty and compact.
Assume that $\Delta$ has a simple zero at every boundary point
of $K$. Then, for every $h\in\R[x]$ strictly positive on $\R$,
there exist $\sigma_0,\sigma_1\in\SOS_{\C}$ such that
\begin{equation}\label{eq:scalar}
 \Delta=\sigma_{0}+h\Delta \sigma_{1}.
\end{equation}
\end{lemma}

\begin{proof}
The polynomial $h\Delta$ defines the same set $K$ and has the same
derivative signs at its boundary as $\Delta$. By
Theorem~\ref{thm:scalar}, $M(\{h\Delta\})$ is saturated. Apply this
to the polynomial $\Delta$, which is nonnegative on $K$.
\end{proof}

We use the following cancellation argument for a strictly positive multiplier.

\begin{lemma}\label{lem:cancel}
Let $h,\Delta\in\R[x]$ satisfy $h>0$ on $\R$ and
\[
\Delta=\sigma_0+h\Delta\sigma_1,
\qquad \sigma_0,\sigma_1\in\SOS_{\C}.
\]
If $F\in\B(\Hh)[x]$ is hermitian and $hF\in\M_{\Hh}(\{\Delta\})_{\le n}$, then
\[
F\in\M_{\Hh}(\{\Delta\})_{\le
\max\{\deg F,n+\deg\sigma_1\}}.
\]
\end{lemma}

\begin{proof}
Write $hF=A^*A+\Delta B^*B$, with $2\deg A\le n$ and $\deg\Delta+2\deg B\le n$. Substituting the scalar identity gives
\begin{equation}\label{eq:cancel-positivity}
h\bigl(F-\Delta\sigma_1B^*B\bigr)
=A^*A+\sigma_0B^*B\succeq0
\quad\text{on }\R.
\end{equation}
Since $h>0$ on $\R$, the hermitian polynomial $G=F-\Delta\sigma_1B^*B$ is positive semidefinite on $\R$. By Proposition~\ref{prop:factorization}, $G=H^*H$ for some $H\in\B(\Hh)[x]$. Thus
\[
F=H^*H+\Delta\sigma_1B^*B \, \in\, \mathcal M_{\mathcal H}(\{\Delta\}).
\]
The degree of the second summand is at most $n+\deg\sigma_1$, while $\deg G\le\max\{\deg F,n+\deg\sigma_1\}$. This gives the stated bound.
\end{proof}

\begin{remark}\label{rem:cancel}
The following remarks explain when nonnegativity of the multiplier $h$ suffices and when strict positivity is required. They should be read alongside Subsection~\ref{ssec:strict-multiplier}.
\begin{enumerate}[leftmargin=*,font=\normalfont,itemsep=3pt]
\item[(a)]
The hypothesis $h>0$ on $\R$ in
Lemma~\ref{lem:scalar-identity} cannot in general be
replaced by $h\ge0$.
Indeed, if $h(c)=0$ at a point where $\Delta(c)<0$,
then the identity
$\Delta=\sigma_0+h\Delta\sigma_1$,
with $\sigma_0,\sigma_1\in\SOS_{\C}$,
would imply $\Delta(c)=\sigma_0(c)\ge0$, a contradiction.
This obstruction can occur even if $h$ is strictly
positive on $K$.
A zero of $h$ at a boundary point $a$ of $K$
also prevents this identity. Since
$\Delta(a)=0$ and $\Delta'(a)\ne0$, the identity
would imply
$\sigma_0(a)=0$ and
$\sigma_0'(a)=\Delta'(a)\ne0$,
which is impossible for a nonnegative polynomial.

\item[(b)]
For the cancellation of $h$ in
\eqref{eq:cancel-positivity}, however, it suffices
to assume that $h\ge0$ on $\R$ and $h\not\equiv0$.
Indeed, this identity implies
$F-\Delta\sigma_1B^*B\succeq0$ wherever $h>0$.
Since a nonzero polynomial has only finitely many
real zeros, the same inequality holds on 
$\R$ by continuity. Thus strict positivity is used
to guarantee the scalar identity in
Lemma~\ref{lem:scalar-identity}, whereas the
cancellation step itself only requires a nonzero
nonnegative multiplier.
\end{enumerate}
\end{remark}

\section{Reduction}\label{sec:reduction}
\subsection{Finite sets}
Theorem~\ref{thmA} admits a direct proof when $K_{S}$ is finite using the Lagrange interpolation polynomials.
\begin{proposition}\label{prop:finite-set}
Let $K=\{\xi_1,\ldots,\xi_r\}\subset\R$, where $r\ge1$ and the
points are distinct. Let $S\subset\R[x]$ be finite with
$K_S=K$ and $M(S)$ saturated.
There is a nonnegative integer $c_S$  such that, for every 
complex Hilbert space $\Hh$, every nonnegative integer $\vd$, and every hermitian
$F\in\B(\Hh)[x]$,
\[
F\succeq0\text{ on }K,\qquad \deg F\le \vd
\quad\Longrightarrow\quad
F\in\M_{\Hh}(S)_{\le \vd+c_{S}}.
\]
Moreover, for the natural description $S_K$ of $K$,
\[
F\succeq0\text{ on }K,\qquad \deg F\le\vd
\quad\Longrightarrow\quad
F\in\mathcal T_{\Hh}(S_K)_{\le\max\{\vd,r-1\}}.
\]
If $r\le3$, the bound $\max\{\vd,r-1\}$ can be
replaced by $\vd$.
\end{proposition}

\begin{proof}
Let $\ell_i$ be the Lagrange interpolation polynomials
\begin{equation}\label{eq:Lagrange}
\ell_i(x)=\prod_{\substack{1\le j\le r\\j\ne i}}\frac{x-\xi_j}{\xi_i-\xi_j}
\qquad(1\le i\le r),
\end{equation}
and set
\[
\varphi(x)=\prod_{i=1}^{r}(x-\xi_i),
\qquad
H(x)=\sum_{i=1}^{r}\ell_i(x)^2F(\xi_i).
\]
Since $F(\xi_i)\succeq0$ for every $i$,
we have $H\in\SOS_{\Hh}$.
Moreover, $H(\xi_i)=F(\xi_i)$, so $F-H$ vanishes at every $\xi_i$ and is therefore divisible by $\varphi$. Consequently,
$ F-H=\varphi G $
for some $G\in\B(\Hh)[x]$. Since $F-H$ is hermitian and $\varphi$ has real coefficients, $G$ is also hermitian.

\medskip

If $G=0$, the assertion is immediate.
Otherwise,
\[
\deg G\le\max\{\vd,2r-2\}-r.
\]
Set $k=\lceil\deg G/2\rceil$.
Since $\deg G\le2k$, we can choose $N>0$
sufficiently large that
\[
G_+(x) =G(x)+N(1+x^2)^kI_{\Hh} \,\succeq\, 0,
\qquad (x \in \R).
\]
The polynomial
\[
G_-= N(1+x^2)^kI_{\Hh}
\]
is also positive semidefinite on $\R$.
By Proposition~\ref{prop:factorization},
both $G_+$ and $G_-$ belong to
$\SOS_{\Hh}$, with
\[
\deg G_\pm\le2k\le\deg G+1.
\]

Now
\[
F=H+\varphi G_++(-\varphi)G_-.
\]
Since $\varphi$ and $-\varphi$ vanish
on $K$, scalar saturation gives
$\varphi,-\varphi\in M(S)$.
Fix scalar certificates for both
polynomials of weighted degree at most $\widetilde c_S$.
Substituting these certificates yields
$ F\in\M_{\Hh}(S)_{\le n}$,
where
\[
n\le
\max\bigl\{
2r-2,\,
\widetilde c_S+\max\{\vd,2r-2\}-r+1
\bigr\}.
\]
Therefore,
$ n\le \vd+c_S$, where 
$ c_S=\max\{2r-2, \widetilde c_S+r-1\}. $

\smallskip

For the additional assertion concerning the natural description $S_K$, write
$F=\sum_{j=0}^{\vd}A_jx^j$ and use the Lagrange
polynomials $\ell_i$ from \eqref{eq:Lagrange}. Let
\[
h_j(x)=x^j-\sum_{i=1}^{r}\xi_i^j\ell_i(x).
\]
Note that $h_j(\xi_i)=0$ for $i=1,\ldots,r$.
Thus $h_j=0$ for $j<r$, and $\deg h_j=j$ for $j\geq r$.
Hence
we can write
\begin{align*}
F(x)
    &=\sum_{j=0}^{\vd}A_j
    \Bigl(h_j(x)+\sum_{i=1}^r \xi_i^j\ell_i(x)\Bigr)\\
    &=\sum_{j=0}^{\vd}A_jh_j(x)
    +\sum_{i=1}^{r}\ell_i(x)
    \Bigl(\sum_{j=0}^{\vd}A_j\xi_i^j\Bigr)\\
    &=\sum_{j=r}^{\vd}A_jh_j(x)
    +\sum_{i=1}^{r}\ell_i(x)F(\xi_i).
\end{align*}
Each $\ell_i$ is nonnegative on $K$ and has degree
at most $r-1$, while $h_j$ vanishes on $K$ and has
degree at most $j$ for $j\ge r$.
The scalar degree-preserving theorem \cite{KM,KMS}
therefore gives natural-preordering certificates
for $\ell_i$ of weighted degree at most $r-1$
and for both $h_j$ and $-h_j$ of weighted degree
at most $j$.

Write $A_j=A_j^+-A_j^-$, where $A_j^+,A_j^-\succeq 0$, and
\[
h_jA_j=h_jA_j^++(-h_j)A_j^-.
\]
Since $F(\xi_i)$ and $A_j^\pm$ are positive semidefinite,
multiplying the scalar certificates by these constant
operators gives
\[
F\in\mathcal T_{\Hh}(S_K)_{\le\max\{\vd,r-1\}}.
\]

If $\vd=0$, then $F$ is a constant positive operator and hence a
hermitian square. For $\vd>0$ and $r\leq2$, we have
$\max\{\vd,r-1\}=\vd$.
Finally, if $r=3$ and $\vd=1$, then $F$ is affine. Setting
$a=\min K$ and $b=\max K$, we obtain
\[
F(x)=\frac{b-x}{b-a}F(a)+\frac{x-a}{b-a}F(b)\in\mathcal T_{\Hh}(S_K)_{\le 1}.
\]
This completes the proof.
\end{proof}

\subsection{Closedness of truncated operator quadratic modules}
The following lemma will be used in the proof of
Theorem~\ref{thm:quadratic-reduction} and Corollary~\ref{cor:degree-preservation}.

\begin{lemma}\label{lem:closed-truncations}
Let $S=\{g_1,\ldots,g_s\}\subseteq\R[x]$ be finite and suppose that $K_S$ has a nonempty interior. Then, for every complex Hilbert space $\Hh$ and every nonnegative integer $n$, the cone $\M_{\Hh}(S)_{\le n}$ is closed under coefficientwise convergence in operator norm.

More generally, if a net
$(F_\lambda)\subseteq\M_{\Hh}(S)_{\le n}$
has uniformly bounded coefficients in operator norm
and converges coefficientwise in the weak operator
topology to $F$, then
$F\in\M_{\Hh}(S)_{\le n}$.
\end{lemma}

\begin{proof}
We prove the more general assertion.
Let $(F_\lambda)\subseteq\M_{\Hh}(S)_{\le n}$ be a net
whose coefficients are uniformly bounded in operator norm
and which converges coefficientwise in the weak operator
topology to $F$.
Choose representations
\begin{align}
\label{eq:closedness-rep}
F_\lambda=\sum_{j=0}^s g_jQ_{j,\lambda},
\qquad
Q_{j,\lambda}\in\SOS_{\Hh},
\qquad g_0=1.
\end{align}
We may omit terms with
$\deg g_j>n$. For the remaining indices, put
\[
r_j=\left\lfloor\frac{n-\deg g_j}{2}\right\rfloor.
\]
Then $\deg Q_{j,\lambda}\le2r_j$, and
$Q_{j,\lambda}\succeq0$ on $\R$.

\smallskip

Choose a nondegenerate compact interval $J\subseteq K_S$
on which all the remaining generators $g_j$ are strictly
positive. 
Since the coefficients of $F_\lambda$ are uniformly
bounded in operator norm, there exists $M>0$ such that
\[
0\preceq F_\lambda(x)\preceq M I_{\Hh}
\qquad(x\in J)
\]
for every $\lambda$.
Each summand in \eqref{eq:closedness-rep} is positive semidefinite on $J$, so
\[
0\preceq Q_{j,\lambda}(x)
 \preceq\frac{M}{g_j(x)}I_{\Hh}
\qquad(x\in J).
\]
Since $g_j$ is bounded below by a positive constant on
$J$, the values of $Q_{j,\lambda}$ are uniformly bounded
in operator norm.

\smallskip

For each $j$, choose $2r_j+1$ distinct points of $J$.
Polynomial interpolation expresses every coefficient of
$Q_{j,\lambda}$ as a fixed linear combination of its
values at these points. Hence all the coefficients of
all the $Q_{j,\lambda}$ are uniformly bounded in
operator norm.

\smallskip

Closed bounded balls of $\B(\Hh)$ are compact in the
weak operator topology. As only finitely many coefficients
occur, a common subnet converges coefficientwise in that
topology to polynomials $Q_j$ with $\deg Q_j\le2r_j$.
Along this subnet, for every $x\in\R$,
$Q_{j,\lambda}(x)$ converges in the weak operator topology
to $Q_j(x)$. 
Thus $Q_j\succeq0$ on $\R$.
By Proposition~\ref{prop:factorization},
each $Q_j$ belongs to $\SOS_{\Hh}$, with a factor of
degree at most $r_j$.
Passing to the limit in the coefficient identities gives
\[
F=\sum_jg_jQ_j,
\]
a certificate of weighted degree at most $n$.
\end{proof}

\subsection{Reduction to quadratic operator polynomials}

\begin{lemma}\label{lem:affine-annihilator}
Let $r\ge1$ and write
$
v(x)=\col(1,x,\ldots,x^r)\in\mathbb C^{r+1}.
$
An affine matrix polynomial
$b(x)=b_0+xb_1\in\mathbb C[x]^{r\times(r+1)}$
satisfies $b(x)v(x)=0$ for all $x\in \mathbb R$ if and only if
\[
b_0=\begin{pmatrix}0&M\end{pmatrix},
\qquad
b_1=\begin{pmatrix}-M&0\end{pmatrix}
\]
for some $M\in\mathbb C^{r\times r}$.
We denote the polynomial determined by $M$ by $b_M(x)$.

Moreover, if $M$ is invertible, then $b_M(x)$ has full row
rank for every $x\in\mathbb R$ and
\begin{equation}\label{eq:b-to-v}
\ker b_M(x)=\{\lambda v(x)\ \colon\ \lambda\in \mathbb C\}
\qquad (x\in \mathbb R).
\end{equation}
\end{lemma}

\begin{proof}
Write $b(x)=b_0+xb_1$, and let
$c_0,\ldots,c_r\in\mathbb C^r$ and
$d_0,\ldots,d_r\in\mathbb C^r$ denote the
column vectors of $b_0$ and $b_1$,
respectively. Thus,
$b_0=\begin{pmatrix}c_0&\cdots&c_r\end{pmatrix}$,
$b_1=\begin{pmatrix}d_0&\cdots&d_r\end{pmatrix}$.
The identity $b(x)v(x)=0$ holds if and only if
\[
c_0+\sum_{j=1}^r(c_j+d_{j-1})x^j
+d_rx^{r+1}=0.
\]
Comparing coefficients, we obtain
$c_0=0$, $d_r=0$ and $c_j=-d_{j-1}$
for $1\le j\le r$.
Thus, with $M=(c_1,\ldots,c_r)$,
we have
$b_0=\begin{pmatrix}0&M\end{pmatrix}$,
$b_1=\begin{pmatrix}-M&0\end{pmatrix}$.

For the moreover part, write
\[
b_M(x)=M
\begin{pmatrix}
-x&1&0&\cdots&0\\
0&-x&1&\ddots&\vdots\\
\vdots&\ddots&\ddots&\ddots&0\\
0&\cdots&0&-x&1
\end{pmatrix}.
\]
Assume that $M$ is invertible.
If $q=\col(q_0,\ldots,q_r)$ belongs to
$\ker b_M(x)$, then
$q_j=xq_{j-1}$ for $1\le j\le r$.
Thus, $q=q_0v(x)$, proving
\eqref{eq:b-to-v}. Since $\dim(\ker b_M(x))=1$, $b_M(x)$ has rank $r$
for every $x\in\mathbb R$.
\end{proof}

\begin{theorem}\label{thm:quadratic-reduction}
Let $S\subseteq\R[x]$ be finite, and suppose that $K_S$ is compact and has a nonempty interior. Assume that there is an integer $N\ge2$ such that for every complex Hilbert space $\mathcal E$ and every hermitian $L\in\B(\mathcal E)[x]$ with $\deg L\le2$,
\[
L\succeq0\text{ on }K_S
\quad\Longrightarrow\quad
L\in\M_{\mathcal E}(S)_{\le N}.
\]
Then, for every complex Hilbert space $\Hh$, every integer $\vd\ge1$, and every hermitian $F\in\B(\Hh)[x]$ with $\deg F\le\vd$,
\[
F\succeq0\text{ on }K_S
\quad\Longrightarrow\quad
F\in\M_{\Hh}(S)_{\le N+2\lceil\vd/2\rceil-2}
\subseteq\M_{\Hh}(S)_{\le\vd+N-1}.
\]
For $\vd=0$, a certificate of degree zero suffices.
\end{theorem}

\begin{proof}
For $\vd=0$, positivity gives the degree-zero
certificate
$F=(F^{1/2})^*F^{1/2}.$
If $1\le\vd\le2$, the assertion follows directly
from the quadratic hypothesis. We therefore
assume $\vd\ge3$ and put
$
r=\left\lceil\frac{\vd}{2}\right\rceil-1\ge1.
$
Write
\[
F(x)=\sum_{k=0}^{2r+2}F_kx^k,
\]
with zero coefficients added if necessary.
Set
\[
\mathcal E=\mathbb C^{r+1}\otimes\Hh
\cong\Hh^{\oplus(r+1)},
\qquad
v(x)=\col(1,x,\ldots,x^r),
\]
and define the polynomial column
\[
V(x)=v(x)\otimes I_{\Hh}:
\Hh\longrightarrow\mathcal E.
\]
Consider the block diagonal hermitian polynomial
\[
A(x)=
\begin{pmatrix}
F_0+xF_1 & 0 & \cdots & 0 & 0\\
0 & F_2+xF_3 & \ddots & \vdots & \vdots\\
\vdots & \ddots & \ddots & 0 & \vdots\\
0 & \cdots & 0 & F_{2r-2}+xF_{2r-1} & 0\\
0 & \cdots & \cdots & 0 & F_{2r}+xF_{2r+1}+x^2F_{2r+2}
\end{pmatrix}.
\]
Then $\deg A\le2$, and direct multiplication gives
\begin{equation}\label{eq:quadratic-compression}
V(x)^*A(x)V(x)=F(x).
\end{equation}

Although $F(x)\succeq0$ on $K_S$, the identity
\eqref{eq:quadratic-compression} only gives
positivity of $A(x)$ on the subspace
$\operatorname{ran}V(x)$, not necessarily on
the whole space $\mathcal E$.
We therefore modify $A$ by a quadratic
positive semidefinite operator polynomial that vanishes
under compression by $V$.

By Lemma~\ref{lem:affine-annihilator}, choosing
$M=I_{\mathbb C^r}$ gives an affine matrix
polynomial
\[
b(x)=b_0+xb_1,\qquad
b_0=\begin{pmatrix}0&I_r\end{pmatrix},
\qquad
b_1=\begin{pmatrix}-I_r&0\end{pmatrix},
\]
such that
\[
b(x)v(x)=0,\qquad
\ker b(x)=\{\lambda v(x)\ \colon \ \lambda \in \mathbb C\}
\]
for every $x\in\mathbb R$.
Set
\[
B(x)=b(x)\otimes I_{\Hh}:
\mathcal E\longrightarrow\Hh^{\oplus r}.
\]
It follows that
\[
B(x)V(x)=0,\qquad
\ker B(x)=\operatorname{ran}V(x).
\]
Since $b(x)$ has full row rank for every
real $x$ and $K_S$ is compact, there exists
$\kappa>0$ such that
\begin{equation}\label{eq:penalty-gap}
\|B(x)\eta\|^2\ge\kappa\|\eta\|^2
\end{equation}
for every $x\in K_S$ and every
$\eta\in(\ker B(x))^\perp$.

Fix $\delta>0$ and set
\[
A_\delta(x)=A(x)+\delta I_{\mathcal E},
\qquad
\Lambda_\delta=
\max_{x\in K_S}\|A_\delta(x)\|.
\]
By \eqref{eq:quadratic-compression},
\[
\begin{aligned}
V(x)^*A_\delta(x)V(x)
&=F(x)+\delta V(x)^*V(x)
\succeq\delta V(x)^*V(x)
\end{aligned}
\]
for every $x\in K_S$.
Thus,
\begin{equation}\label{eq:positive-on-kernel}
\langle A_\delta(x)\xi_0,\xi_0\rangle
\ge\delta\|\xi_0\|^2
\end{equation}
for every $\xi_0\in\ker B(x)$.

\smallskip

We now show that
if we add 
a sufficiently large
multiple of $B(x)^*B(x)$ to $A_\delta(x)$, then $A_\delta(x)$
becomes positive semidefinite on the whole space
$\mathcal E$.
Fix $x\in K_S$ and let $\xi\in\mathcal E$
be arbitrary. We write
\[
\xi=\xi_0+\xi_1,\qquad
\xi_0\in\ker B(x),\quad
\xi_1\in(\ker B(x))^\perp.
\]
Using \eqref{eq:positive-on-kernel} and the
definition of $\Lambda_\delta$, we obtain
\begin{align*}
\langle A_\delta(x)\xi,\xi\rangle
&\ge
\delta\|\xi_0\|^2
-2\Lambda_\delta\|\xi_0\|\|\xi_1\|
-\Lambda_\delta\|\xi_1\|^2
\ge
\frac{\delta}{2}\|\xi_0\|^2
-\left(
\Lambda_\delta+
\frac{2\Lambda_\delta^2}{\delta}
\right)\|\xi_1\|^2.
\end{align*}
Choose
\[
t_\delta\ge
\kappa^{-1}\left(
\Lambda_\delta+
\frac{2\Lambda_\delta^2}{\delta}
\right)
\]
and define
\[
L_\delta(x)
=A_\delta(x)+t_\delta B(x)^*B(x).
\]
Since $B(x)\xi_0=0$, the estimate
\eqref{eq:penalty-gap} yields
\begin{align*}
\langle L_\delta(x)\xi,\xi\rangle
&=
\langle A_\delta(x)\xi,\xi\rangle
+t_\delta\|B(x)\xi_1\|^2\\
&\ge
\frac{\delta}{2}\|\xi_0\|^2
+
\left(
t_\delta\kappa-\Lambda_\delta
-\frac{2\Lambda_\delta^2}{\delta}
\right)\|\xi_1\|^2\\
&\ge0.
\end{align*}
As $\xi\in\mathcal E$ and $x\in K_S$
were arbitrary, we have
$
L_\delta(x)\succeq0$
for $x\in K_S.$
Moreover, $L_\delta$ is a hermitian operator
polynomial of degree at most two. By the hypothesis,
\[
L_\delta\in\M_{\mathcal E}(S)_{\le N}.
\]
Thus $L_\delta$ is of the form
\[
L_\delta
=\sum_{j=0}^s
g_jH_{j,\delta}^*H_{j,\delta},
\qquad g_0=1,
\]
where $H_{j,\delta}\in\B(\mathcal E)[x]$
and
$
\deg(g_jH_{j,\delta}^*H_{j,\delta})
\le N.
$
Since $B(x)V(x)=0$, we have
\begin{align*}
V(x)^*L_\delta(x)V(x)
&=V(x)^*A_\delta(x)V(x)
=F(x)+\delta V(x)^*V(x)
=F(x)+\delta
\left(\sum_{k=0}^r x^{2k}\right)I_{\Hh}.
\end{align*}
On the other hand,
\[
V(x)^*L_\delta(x)V(x)
=\sum_{j=0}^s g_j(x)Q_{j,\delta}(x),
\qquad
\text{where}\quad
Q_{j,\delta}(x)
=V(x)^*H_{j,\delta}(x)^*
H_{j,\delta}(x)V(x).
\]
Each $Q_{j,\delta}\in \B(\mathcal H)[x]$ is positive
semidefinite on $\mathbb R$, and
$
\deg(g_jQ_{j,\delta})\le N+2r.
$
By Proposition~\ref{prop:factorization},
each $Q_{j,\delta}$ belongs to
$\SOS_{\Hh}$, without increasing its
degree. Consequently,
\[
F+\delta
\left(\sum_{k=0}^r x^{2k}\right)I_{\Hh}
\in\M_{\Hh}(S)_{\le N+2r}.
\]

As $\delta\downarrow0$, the left-hand side
converges coefficientwise in operator norm
to $F$. By Lemma~\ref{lem:closed-truncations},
the truncated operator quadratic module
$\M_{\Hh}(S)_{\le N+2r}$ is closed.
Therefore,
\[
F\in\M_{\Hh}(S)_{\le N+2r}
=\M_{\Hh}(S)_{
\le N+2\lceil\vd/2\rceil-2}.
\]
Finally,
$N+2\left\lceil\frac{\vd}{2}\right\rceil-2
\le\vd+N-1,$
which proves the assertion.
\end{proof}

\begin{remark}
\label{rem:scope-quadratic-reduction}
The algebraic construction in the proof of
Theorem~\ref{thm:quadratic-reduction} extends to several
variables. One replaces the univariate
monomial column by the column of all monomials of
total degree at most $r$, together with affine
relations whose common kernel is spanned by this
column.
On a compact basic closed set, provided that the
relevant truncated operator quadratic module is
closed, a uniform degree bound $N$ for certificates
of all positive semidefinite quadratic operator
polynomials would therefore imply the bound
$\vd+N-1$, $\vd\ge1$,
for polynomials of total degree at most $\vd$.
The bound $N$ must be independent of both the
quadratic polynomial and its coefficient Hilbert
space, since the reduction enlarges that space.
\end{remark}

\section{Quadratic certificates on finite unions of intervals}\label{sec:intervals}
Throughout this section, fix $\ell\in \mathbb N$,
\begin{equation}\label{eq:endpoint-data}
 K \,=\, \bigcup_{j=1}^{\ell}[a_j,b_j],\qquad
 a_1<b_1<\cdots<a_\ell<b_\ell,
\end{equation}
and define
\begin{equation}\label{eq:P-Q-r-Delta}
 \begin{aligned}
 P(x)& \,=\, \prod_{j=1}^{\ell}(x-a_j),\qquad&
 Q(x)& \,=\, \prod_{j=1}^{\ell}(x-b_j),\\
 R(x)& \,=\, Q(x)-P(x),\qquad&
 \Delta(x)& \,=\, -P(x)Q(x).
 \end{aligned}
\end{equation}
The set $K$ is $\{\Delta\ge0\}$, with simple boundary zeros, so
$M(\{\Delta\})$ is saturated by Theorem~\ref{thm:scalar}.

\subsection{A symmetric matrix pencil}

Our aim is to reduce polynomial positivity on $K$
to positivity on $[0,1]$, where the operator
Fej\'er--Riesz theorem provides certificates
with controlled degree. A natural first choice
is the diagonal pencil
\[
\widetilde T(y)
=(1-y)\diag(b_1,\ldots,b_\ell)
+y\diag(a_1,\ldots,a_\ell).
\]
Its eigenvalues cover the intervals
$[a_j,b_j]$ as $y$ runs over $[0,1]$.
However, recovering $x$ from the $j$th
eigenvalue requires the substitution
$ y=(b_j-x)/(b_j-a_j),$
which depends on $j$. Thus, this choice
does not provide a common rational parameter
for all the intervals.

\medskip

We therefore seek a linear pencil $T(y)$ whose
characteristic polynomial is affine in $y$,
with $Q$ and $P$ as its characteristic
polynomials at $y=0$ and $y=1$, respectively.
Thus, we require
\[
\det(xI-T(y))
=(1-y)Q(x)+yP(x)
=Q(x)-y R(x).
\]

A rank-one perturbation has the required
affine dependence on $y$. We therefore look
for $v\in\R^\ell$ such that
\begin{equation}\label{def:T} 
T(y)=D-yvv^{\mathsf T},
\qquad D=\diag(b_1,\ldots,b_\ell),
\end{equation}
has eigenvalues $a_1,\ldots,a_\ell$ at $y=1$.
The following lemma shows that such a vector exists.

\begin{lemma}\label{lem:pencil}
There is a column vector $v = \col (v_{1}, \ldots, v_{\ell})\in\R^\ell$, with every $v_j\ne0$, such that the one-parameter rank-one perturbation $T(y)$ defined in \eqref{def:T}
satisfies
\begin{equation}\label{eq:characteristic}
 \det(xI-T(y))=Q(x)-y R(x),\qquad
 a_j\le\lambda_j(T(y))\le b_j\quad(0\le y\le1),
\end{equation}
where eigenvalues are ordered increasingly.
\end{lemma}

\begin{proof}
Applying \cite[Theorem 4.3.26]{HJ} to the reversed ordering
$-b_\ell<-a_\ell<\cdots<-b_1<-a_1$, we obtain
$v={\rm col}(v_1,\ldots,v_\ell)\in\R^\ell$ such that
$D-vv^{\mathsf T}$ has eigenvalues
$\{a_1,\ldots,a_\ell\}$.
Each $v_j$ is nonzero, since otherwise $b_j$
would remain an eigenvalue, contradicting
the prescribed spectrum.

\medskip

Recall the rank-one determinant identity
\cite[Equation~(0.8.5.11)]{HJ},
\begin{equation}\label{eq:rank-one determinant formula}
 \det(A+uv^{\mathsf T})
 =\det A+v^{\mathsf T}\operatorname{adj}(A)u.
\end{equation}
Taking $A=xI-D$ and $u=yv \in \R^{\ell}$, we obtain
\[
 \det(xI-T(y))
 =Q(x)+y\,v^{\mathsf T}\adj(xI-D)v.
\]
At $y=1$, the left-hand side equals $P(x)$, so
\begin{equation}
\label{eq:v^Tu}
v^{\mathsf T}\adj(xI-D)v=P(x)-Q(x)=-R(x).
\end{equation}
Consequently,
$ \det(xI-T(y))=Q(x)-yR(x), $
which proves the first identity.
The inequalities follow
from $T(1)\preceq T(y)\preceq T(0)$.
\end{proof}

Fix the pencil $T(y)$ as in \eqref{def:T}. The eigenvector equation
$T(y)u(x) = x u(x)$ requires $(xI-D)u(x)$ to be a scalar multiple of $v$.
Away from the $b_j$, this suggests the direction $(xI-D)^{-1}v$.
Clearing the denominator $Q(x)$ leads to the polynomial vector
\begin{equation}\label{eq:u}
 u(x)=\adj(xI-D)v.
\end{equation}
The following lemma shows that $u(x)$ never vanishes on $\R$
and establishes an identity expressing the parameter $y$
as a rational function of $x$.

\begin{lemma}\label{lem:recovery-vector}
The polynomial vector $u(x)$ defined in \eqref{eq:u} has no real zero and satisfies
\begin{equation}\label{eq:recovery-identity}
 v^{\mathsf T}u(x)=-R(x),\qquad
 (T(y)-xI)u(x)=(yR(x)-Q(x))v.
\end{equation}
The degree of $u$ is $\ell-1$, and its leading coefficient is $v$.
Consequently, the scalar polynomial 
\begin{equation} 
\label{eq:definition-w}
    w=u^*u
\end{equation}
satisfies
\begin{equation*}\label{eq:w}
 w>0\text{ on }\R,\qquad \deg w=2\ell-2.
\end{equation*}
\end{lemma}

\begin{proof}
The first equality in \eqref{eq:recovery-identity} follows directly from \eqref{eq:v^Tu}
and \eqref{eq:u}.
The adjugate identity $(xI-D)\adj(xI-D)= \det (xI-D) I  = Q(x)I$ also gives
\[
(D-xI)u(x)= (D-xI) \adj (xI-D)v = -Q(x)v.
\]
Therefore,
\begin{align*}
 (T(y)-xI)u(x)
 & \,=\,
(D-yvv^{\mathsf T} - xI) u(x)\\
& \, =\, 
 (D-xI)u(x)-yv(v^{\mathsf T}u(x)) \\
& \,=\, (yR(x)-Q(x))v.
\end{align*}
This proves \eqref{eq:recovery-identity}.

\medskip

Since $xI-D$ is diagonal, the entries of $u(x)$ are
\begin{equation}\label{eq:comp-formula}
 u_j(x)=v_j\prod_{i\ne j}(x-b_i),\qquad j=1,\ldots,\ell.
\end{equation}
If $x\notin\{b_1,\ldots,b_\ell\}$, each entry of $u(x)$ is nonzero. On the other hand, at $x=b_j$, only the $j$th component can be nonzero, and
\[
u(b_j)=v_j\prod_{i\ne j}(b_j-b_i)e_j
=Q'(b_j)v_je_j\ne0,
\]
where $e_j$ is the $j$th standard basis vector and $Q'$ denotes the derivative of $Q$. The last inequality follows from the fact that the $b_i$ are distinct and $v_j\ne0$. Hence $u(x)\ne0$ for every $x\in\mathbb R$.

\medskip

The formulas \eqref{eq:comp-formula}
show that $u$ has degree $\ell-1$ and leading coefficient $v$.
Since $u(x)\ne0$ for every $x\in\R$, the polynomial
$w(x)=\sum_{j=1}^{\ell}u_j(x)^2$ is strictly positive on $\R$.
It has degree $2\ell-2$ and leading coefficient
$v^{\mathsf T}v>0$.
\end{proof}
We now express the parameter $y$ in terms of $x$. Motivated by \eqref{eq:recovery-identity}, we set
\begin{equation}\label{eq:parameter}
y=t(x):=\frac{Q(x)}{R(x)},
\qquad T(t(x))u(x)=xu(x).
\end{equation}
This gives a single rational substitution for all the interval components.

\begin{lemma}\label{lem:parameter}
The polynomial $R$ defined in \eqref{eq:P-Q-r-Delta} has degree
$\ell-1$ and exactly one simple zero in each gap
$(b_j,a_{j+1})$, and no other zeros. In particular, $t(x)$ is
well-defined on $K$ and satisfies
\begin{equation}\label{eq:weight-pullback}
t(x)(1-t(x))=\frac{\Delta(x)}{R(x)^2},
\qquad t(K)\subseteq[0,1],
\qquad t(a_j)=1,\quad t(b_j)=0.
\end{equation}
\end{lemma}
\begin{proof}
Recall that $R(x) = Q(x) - P(x).$ Since both $P$ and $Q$ are monic polynomials of degree $\ell$, the polynomial $R$ has degree at most $\ell -1$. The coefficient of $x^{\ell-1}$ in $R$ is
$\sum_j a_j-\sum_j b_j<0$. Thus the degree of $R$ is exactly $\ell-1.$  
 Since
$R(b_j)=-P(b_j)$ and $R(a_{j+1})=Q(a_{j+1})$ have opposite signs, each gap $(b_{j}, a_{j+1})$ has at least one zero of $R.$ 
The resulting $\ell-1$ distinct zeros exhaust the degree.
The remaining assertions follow from
$Q(R-Q)=-PQ=\Delta$ and the endpoint values.
\end{proof}

For $F(x)=\sum_k A_{k}x^k $, put
\[
 F(T(y)) = \sum_k A_{k} \otimes T(y)^k,\qquad
 \mathbf u(x)= I_{\Hh} \otimes u(x),
\]
where $I_{\Hh}$ denotes the identity operator on $\Hh$.

\begin{lemma}\label{lem:compression}
Let $F\in\B(\Hh)[x]$. Then
\begin{equation}\label{eq:compression}
F(T(t(x)))\mathbf u(x)=\mathbf u(x) F(x),
\qquad
\mathbf u(x)^*F(T(t(x)))\mathbf u(x)=w(x)F(x)
\end{equation}
hold as rational identities. Moreover, if $F$ is hermitian and
$F\succeq0$ on $K$, then
\[
F(T(y))\succeq0 \qquad \text{for all }y\in[0,1].
\]
\end{lemma}

\begin{proof}
Fix $0 \leq y \leq 1$. Since the matrix $T(y)$ is real symmetric, there exists an orthogonal matrix $U_y\in\R^{\ell\times\ell}$ such that
$ U_y^{\mathsf T}T(y)U_y
 =\diag(\lambda_1,\ldots,\lambda_\ell). $
By Lemma~\ref{lem:pencil},
$\lambda_j\in[a_j,b_j]\subseteq K$, so
\[
 (I_{\Hh} \otimes U_y^{\mathsf T}) F(T(y))
 (I_{\Hh} \otimes U_y)
\,=\,
\bigoplus_{j=1}^{\ell}F(\lambda_j)
\,\succeq \, 0.
\]
Hence $F(T(y))\succeq 0$.

\medskip

Write $F(x)=\sum_{k=0}^{\vd} A_{k} x^k $.
For $R(x)\ne0$, \eqref{eq:parameter} gives
\[
 (T(t(x)))^k u(x)=x^k u(x)\qquad(k\ge0).
\]
Thus,
 \begin{align*}
 F(T(t(x)))\mathbf u(x)
 & \,=\, \sum_{k=0}^{\vd}
   A_{k} \otimes (T(t(x)))^{k} (I_{\Hh} \otimes u(x))\\
 & \,=\, \sum_{k=0}^{\vd} A_{k} \otimes x^k u(x)
 =\mathbf u(x) F(x) ,
 \end{align*}
and
\[
 \mathbf u(x)^* F(T(t(x))) \mathbf u(x)
 \,=\, 
 \mathbf u(x)^*\mathbf u(x) F(x)
 \,=\,
 (u(x)^{\mathsf T} u(x))F(x)
 \,=\, w(x)F(x).
\]
These are identities of rational operator functions.
\end{proof}

\subsection{A strictly positive multiplier}
\label{ssec:strict-multiplier}
Let $F\in\B(\Hh)[x]$ be hermitian, with $\deg F\le 2$
and $F\succeq0$ on $K$. The lifted polynomial
$L(y)\coloneqq F(T(y))$ is quadratic and positive
semidefinite on $[0,1]$. By Proposition~\ref{prop:factorization},
it admits an interval certificate. Substituting
$y=t(x)=Q(x)/R(x)$, clearing denominators, and compressing
by $\mathbf u(x)$ gives a certificate for $R^2wF$ by
\eqref{eq:compression}. 
The multiplier $R^2w$ vanishes at
the zeros of $R$, which lie in the gaps between the interval components of $K$. Lemma~\ref{lem:cancel}
(see also Remark~\ref{rem:cancel}) therefore does not apply directly.
To obtain a strictly positive multiplier, we modify the lifted polynomial by a polynomial congruence 
\[
G_{C}(y)=(I_{\Hh}\otimes C(y)^*)F(T(y))
     (I_{\Hh}\otimes C(y)),
\]
which preserves positivity on $[0,1]$.
As in \eqref{eq:compression}, we want to recover a scalar
multiple of $F(x)$ by compressing $G_{C}(t(x))$.
Recall from \eqref{eq:parameter} that
$T(t(x))u(x)=xu(x)$. It is this eigenvector relation
that allows us to recover $F(x)$ from $F(T(t(x)))$.
To retain this property after the congruence, we require
\begin{equation}\label{eq:recovery es}
C(t(x))u(x) = \beta(x)u(x)
\end{equation}
for some scalar function $\beta$. We then have
\[
\mathbf u(x)^*G_{C}(t(x))\mathbf u(x)
=|\beta(x)|^2w(x)F(x).
\]
Thus, we look for a polynomial matrix $C(y)$ that
commutes with $T(y)$, since this ensures that
the eigenspaces of $T(y)$ are preserved.

\medskip

We also want the modified polynomial $G_{C}$ to have degree at most two.
The reason is that increasing the degree of $G_{C}$ may
introduce additional powers of $R$ when we clear
denominators, leaving zeros in the resulting multiplier.
To see this, first consider the original lifted polynomial
$L(y)=F(T(y))$. Since $L$ is quadratic, its interval
factorization has the form
\[
L(y)=H(y)^*H(y)+y(1-y)A^*A,
\qquad H(y)=H_0+yH_1,
\]
where $A$ is constant.
Substituting $y=t(x)=Q(x)/R(x)$
and using $t(1-t)=\Delta/R^2$, we obtain
\[
R^2wF
=R^2\mathbf u^*L \mathbf u
=\bigl((R H_0 + Q H_1)\mathbf u\bigr)^*
\bigl((R H_0+ Q H_1)\mathbf u\bigr)
+\Delta(A\mathbf u)^*(A\mathbf u).
\]
Thus, $R^2$ is enough to clear the denominators.
Now suppose that $C(y)=C_0+yC_1$ is affine.
Even then, the congruence $G_{C} = C^*LC$ may have
degree four. Its interval factorization takes the form
\[
G_{C}(y)=M(y)^*M(y)+y(1-y)N(y)^*N(y),
\]
where
\[
M(y)=M_0+yM_1+y^2M_2,
\qquad N(y)=N_0+yN_1.
\]
This time, clearing the denominators generally
requires multiplication by $R^4$. Indeed, setting
\[
\widetilde M
=R^2M =R^2M_0+RQM_1+Q^2M_2, \qquad
\widetilde N
=RN =RN_0+QN_1,
\]
we obtain
\[
R^4\mathbf u^*G_{C}(t)\mathbf u
=(\widetilde M\mathbf u)^*
 (\widetilde M\mathbf u)
+\Delta(\widetilde N\mathbf u)^*
       (\widetilde N\mathbf u).
\]
Suppose that $C$ satisfies \eqref{eq:recovery es}.
The multiplier in the above certificate is then
\[
R^4|\beta|^2w=R^2|R\beta|^2w.
\]
Since $C$ is affine,
\[
R(x)C(t(x))=R(x)C_0+Q(x)C_1
\]
is a polynomial. Moreover,
\[
(RC_0+QC_1)u=(R\beta)u.
\]
As $u$ never vanishes on $\R$, the scalar function
$R\beta$ has no poles at the real zeros of $R$.
Thus, the multiplier $R^2|R\beta|^2w$ still vanishes
at the zeros of $R$.

\medskip

We therefore seek a polynomial matrix $C$ satisfying
\eqref{eq:recovery es} for some scalar function $\beta$,
such that $\deg G\le 2$ and $R\beta$ is a polynomial
with no real zeros. The following lemma characterizes all such matrix polynomials.

\begin{lemma} \label{lem:degree-preserving-congruences}
Let $\ell\ge 2$. For $C \in \mathbb C[y]^{\ell \times \ell}$, the following
are equivalent:
\begin{enumerate}[label=\textup{(\roman*)},itemsep=3pt]
\item For every $x\in K$, the vector $C(t(x))u(x)$ is a scalar
multiple of $u(x)$, and for every quadratic hermitian
operator polynomial $F$, the polynomial
\[
G_C(y) = (I_{\mathcal H}\otimes C(y)^*) F(T(y))(I_{\mathcal H}\otimes C(y))
\]
has degree at most two.

\item
$C$ commutes with $T$, and both $C$ and $TC$ are affine
in $y$.

\item
There is a polynomial $\vq\in\mathbb C[x]$,
with $\deg \vq\le \ell-1$, such that
$ C(y) = \vq(T(y))R(T(y))^{-1}$.
\end{enumerate}
\end{lemma}

\begin{proof}
We first show that \textup{(i)} $\Longrightarrow$ \textup{(ii)}.
Taking $F=1$ and $F=x^2$, we obtain
\[
\deg(C^*C)\le2,\qquad
\deg((TC)^*(TC))\le2.
\]
Since $\deg(M^*M)=2\deg M$ for every nonzero
matrix polynomial $M$, it follows that both $C$ and
$TC$ have degree at most one.

\smallskip

Next, fix $y\in(0,1)$. The eigenvalues
$\lambda_1(y),\ldots,\lambda_\ell(y)$ of $T(y)$
are distinct and belong to $K$. Since
$t(\lambda_j(y))=y$, the nonzero vectors
$u(\lambda_j(y))$ form an eigenbasis of $T(y)$.
By assumption, $C(y)$ preserves each of these
eigenspaces, so $C(y)$ commutes with $T(y)$.
As this holds for every $y\in(0,1)$ and both
$C$ and $T$ are matrix polynomials, we conclude
that $CT=TC$ identically.

\smallskip

We next establish \textup{(iii)} $\Longrightarrow$ \textup{(ii)}.
Let $\xi_1,\ldots,\xi_{\ell-1}$ be the simple real
zeros of $R$. By \eqref{eq:characteristic},
\[
\det(\xi_jI-T(y))
=Q(\xi_j)-yR(\xi_j)=Q(\xi_j)\neq0.
\]
Thus $T(y)-\xi_jI$ is invertible for every $y$.
Set $A_j=D-\xi_jI$, which is invertible since
$Q(\xi_j)\neq0$. By identity \eqref{eq:rank-one determinant formula},
\[
\det(T(y)-\xi_jI)
= \det(D - yv v^{\mathsf T} - \xi_{j} I) = \det(A_j)(1-yv^{\mathsf T} A_j^{-1}v).
\]
Since $\det(T(y)-\xi_jI)$ is independent of $y$,
we must have
\[
v^{\mathsf T} A_j^{-1}v=0.
\]
The Sherman--Morrison formula therefore gives
\[
(T(y)-\xi_jI)^{-1}
=(A_j-y v v^{\mathsf T})^{-1}
=A_j^{-1}+y A_j^{-1} v v^{\mathsf T} A_j^{-1}.
\]
Thus $ E_{j}(y) \coloneqq (T(y)-\xi_jI)^{-1}$ is affine in $y$. 
The identity
\[
T(y) E_{j}(y) = I+\xi_j E_{j}(y)
\]
shows that $T E_{j} $ is affine as well.

\smallskip

For $\deg \vq\le\ell-1$,
\[
\frac{\vq(x)}{R(x)}
=c+\sum_{j=1}^{\ell-1}
\frac{\vq(\xi_j)}{R'(\xi_j)(x-\xi_j)}.
\]
Substituting $x = T(y)$ yields
\[
C_{\vq}(y)
=cI+\sum_{j=1}^{\ell-1}
\frac{\vq(\xi_j)}{R'(\xi_j)}E_j(y).
\]
Hence $C_{\vq}$ and $T C_{\vq}$ are affine, and $C_{\vq}$
commutes with $T$. Since $C=C_{\vq}$ by assumption, this proves
\textup{(iii)} $\Longrightarrow$ \textup{(ii)}.

\smallskip

Now, we prove \textup{(ii)} $\Longrightarrow$ \textup{(iii)}. Write
$C(y) = B_0 + y B_1$. Since $CT = TC$, we have
$B_0 D = D B_0$. As the diagonal entries of $D$ are
distinct, $B_0$ must be diagonal, say
\[
B_0=\operatorname{diag}(\alpha_1,\ldots,\alpha_\ell).
\]
Choose $\vq\in\mathbb C[x]_{\le\ell-1}$
interpolating
\[
\vq(b_j)=R(b_j)\alpha_j,
\qquad j=1,\ldots,\ell.
\]
Then $C_{\vq}(0)= \vq (D)R(D)^{-1}=B_0$.
The difference $H = C-C_{\vq}$ is therefore of the
form $H(y)=yH_1$. Since $H$ commutes with $T$,
comparison of coefficients in $HT = TH$ gives
$H_{1} D = D H_{1}$, so $H_1$ is diagonal. Furthermore,
$TH$ is affine, whereas
\[
T(y)H(y) = y D H_1-y^2vv^{\mathsf T} H_1.
\]
Thus $vv^{\mathsf T} H_1 = 0$. As every entry of $v$
is nonzero, $H_1=0$. Hence $C = C_{\vq}$, proving
\textup{(ii)} $\Longrightarrow$ \textup{(iii)}.

Finally, \textup{(ii)} $\Longrightarrow$ \textup{(i)}. 
Commutativity preserves the one-dimensional
eigenspaces of $T(t(x))$, while the expressions
$ C^*C, C^*(TC)$ and $ (TC)^*(TC) $ have degree at most two. This proves \textup{(ii)} $\Longrightarrow$ \textup{(i)}. 
\end{proof}

\begin{corollary}\label{cor:degree-preserving-congruences}
For every $\vq\in\C[x]$ with
$\deg\vq\le\ell-1$, the matrix
\begin{equation}\label{eq:C_q}
    C_{\vq}(y):=\vq(T(y))R(T(y))^{-1}
\end{equation}
is affine. Moreover, for every hermitian
$F\in\B(\Hh)[x]$ with $\deg F\le2$,
the polynomial
\begin{equation}\label{eq:lift}
G_{C_{\vq}}(y)
=(I_{\Hh}\otimes C_{\vq}(y)^*)
 F(T(y))(I_{\Hh}\otimes C_{\vq}(y))
\end{equation}
has degree at most two and satisfies
\begin{equation}\label{eq:adjugate-eigenvector}
R(x)C_{\vq}(t(x))u(x)=\vq(x)u(x),
\end{equation}
and
\[
R(x)^2\mathbf u(x)^*
G_{C_{\vq}}(t(x))\mathbf u(x)
=w(x)|\vq(x)|^2F(x).
\]
These identities extend polynomially
across the zeros of $R$.
In particular, the multiplier
$w|\vq|^2$ is strictly positive
on $\R$ if and only if $\vq$ has no
real zeros.
\end{corollary}

\begin{proof}
If $\ell=1$, then $R=a_1-b_1\ne0$ and $\vq$ are constant,
so the assertions follow directly. We may therefore assume
that $\ell\ge2$.
By Lemma~\ref{lem:degree-preserving-congruences},
$C_{\vq}$ is affine in $y$, and $G_{C_{\vq}}$
has degree at most two. For $C=C_{\vq}$, the identity \eqref{eq:parameter} gives
\[
C(t(x))u(x)=\frac{\vq(x)}{R(x)}u(x)
\]
whenever $R(x)\neq0$. Therefore,
\[
R(x)^2\mathbf u^*(x)G_C(t(x))\mathbf u(x)
=|\vq(x)|^2\mathbf u^*(x)F(T(t(x)))\mathbf u(x)
=w(x)|\vq(x)|^2F(x).
\]
Both sides extend polynomially across the zeros
of $R$. Since $w>0$ on $\mathbb R$, the multiplier
$w|\vq|^2$ is strictly positive precisely when
$\vq$ has no real zeros.
\end{proof}

\begin{remark}
An upper bound for the constant $c_S$ in Theorem~\ref{thm:main}
depends on the choice of $\vq$. It would be interesting
to study this dependence and optimize over the family
of multipliers $ w |\vq|^{2}$. We do not pursue this question
here, as any polynomial $\vq$ of degree at most $\ell-1$
with no real zeros suffices for our purposes. Henceforth,
we take $\vq=1$ and write $C=C_{\vq}$ and $G=G_{C_{\vq}}$
for this choice.
\end{remark}

We apply interval factorization to $G$, substitute $y=t(x)$, and compress by $\mathbf u$.

\begin{lemma}\label{lem:pullback}
Let $F\in\B(\Hh)[x]$ be hermitian, with $\deg F\le2$
and $F\succeq0$ on $K$. Let $w$ be defined by \eqref{eq:definition-w}. There exist polynomial columns
\begin{equation}\label{eq:UV-columns}
U=\col(U_1,\ldots,U_\ell),\qquad
V=\col(V_1,\ldots,V_\ell),
\end{equation}
where $U_i,V_i\in\B(\Hh)[x]$ for $i=1,\ldots,\ell$,
with $\deg U\le2\ell-1$ and $\deg V\le\ell-1$,
such that
\begin{equation}\label{eq:raw-certificate}
wF=U^*U+\Delta V^*V\in\M_{\Hh}(\{\Delta\})_{\le4\ell-2}.
\end{equation}
\end{lemma}

\begin{proof}
By Corollary~\ref{cor:degree-preserving-congruences} and \ref{item:interval-factorization} of Proposition~\ref{prop:factorization}, there exist an affine operator polynomial
$G_1\in\B(\Hh^{\oplus\ell})[y]$ and an operator
$H_1\in\B(\Hh^{\oplus\ell})$ such that
\[
G(y)=G_1(y)^*G_1(y)+y(1-y)H_1^*H_1.
\]
Recall that $\mathbf u(x)=I_{\Hh}\otimes u(x)$, where
$u(x)=\adj(xI-D)v$ is the polynomial vector from \eqref{eq:u}. Write $G_1(y)=A+yB$, where
$A,B\in\B(\Hh^{\oplus\ell})$, and set
\[
U(x)=R(x)G_1(t(x))\mathbf u(x)
     =\bigl(R(x)A+Q(x)B\bigr)\mathbf u(x),\qquad
V(x)=H_1\mathbf u(x).
\]
Thus $U$ and $V$ are polynomial columns of the form
\eqref{eq:UV-columns}. Since
$\deg Q=\ell$ and $\deg R=\deg\mathbf u=\ell-1$,
we have $\deg U\le2\ell-1$ and $\deg V\le\ell-1$.

\smallskip

For $x\in\R$ with $R(x)\ne0$, the identity
\eqref{eq:weight-pullback} gives
\[
U(x)^*U(x)+\Delta(x)V(x)^*V(x)
=R(x)^2\mathbf u(x)^*G(t(x))\mathbf u(x).
\]
By \eqref{eq:adjugate-eigenvector},
\[
R(x)(I_{\Hh}\otimes C(t(x)))\mathbf u(x)
=\mathbf u(x).
\]
Hence, using \eqref{eq:lift} and \eqref{eq:compression},
\begin{align*}
R(x)^2\mathbf u(x)^*G(t(x))\mathbf u(x)
&=
  \mathbf u(x)^*F(T(t(x)))\mathbf u(x)
=w(x)F(x).
\end{align*}
Since both sides are polynomial expressions, the identity
\eqref{eq:raw-certificate} holds on all of $\R$.

\smallskip

The degree bounds on $U$ and $V$ yield
$\deg(U^*U)\le4\ell-2$,
$\deg(\Delta V^*V)\le4\ell-2$.
Moreover, $U^*U$ and $V^*V$ are positive semidefinite
on $\R$. By \ref{item:global-factorization} of
Proposition~\ref{prop:factorization}, both belong to
$\SOS_{\Hh}$, with factors of degrees at most
$2\ell-1$ and $\ell-1$, respectively.
Thus,
$
wF\in\M_{\Hh}(\{\Delta\})_{\le4\ell-2}.
$
\end{proof}

\subsection{Quadratic degree bounds}
The polynomial $w$, defined by \eqref{eq:definition-w},
is strictly positive by Lemma~\ref{lem:recovery-vector}.
By Lemma~\ref{lem:scalar-identity}, there are $\sigma_0,\sigma_1\in\SOS_{\C}$ such that
\begin{equation}\label{eq:scalar2}
\Delta=\sigma_0+w\Delta\sigma_1.
\end{equation}
We fix such a pair for the remainder of this section. Note that $\sigma_1\ne0$, since otherwise $\Delta=\sigma_0$ would be nonnegative on all of $\R$.

\begin{theorem}\label{thm:intervals}
Let $K$ and $\ell$ be as in \eqref{eq:endpoint-data}, and let $\Delta$ and $w$ be defined by \eqref{eq:P-Q-r-Delta} and \eqref{eq:definition-w}, respectively. Fix $\sigma_0,\sigma_1\in\SOS_{\C}$ satisfying \eqref{eq:scalar2}. Then, for every complex Hilbert space $\Hh$ and every hermitian $F\in\B(\Hh)[x]$ of degree at most $2$,
\[
F\succeq0\text{ on }K
\quad\Longleftrightarrow\quad
F\in\M_{\Hh}(\{\Delta\})_{\le4\ell-2+\deg\sigma_1}.
\]
\end{theorem}

\begin{proof}
If $F\succeq0$ on $K$, Lemma~\ref{lem:pullback} gives $wF\in\M_{\Hh}(\{\Delta\})_{\le4\ell-2}$. By Lemma~\ref{lem:cancel},
\[
F\in\M_{\Hh}(\{\Delta\})_{\le
\max\{\deg F,4\ell-2+\deg\sigma_1\}}
=\M_{\Hh}(\{\Delta\})_{\le4\ell-2+\deg\sigma_1}.
\]
The converse follows by evaluating the certificate on $K=\{\Delta\ge0\}$.
\end{proof}

\begin{corollary}\label{cor:interval-transfer}
Let $K$ and $\ell$ be as in \eqref{eq:endpoint-data},
and let $\Delta$ and $w$ be defined by
\eqref{eq:P-Q-r-Delta} and \eqref{eq:definition-w}, respectively.
Fix $\sigma_0,\sigma_1\in\SOS_{\C}$ satisfying
\eqref{eq:scalar2}, and let $S$ be a saturated scalar
description of $K$. Then there exists a nonnegative
integer $\eta_S$, depending only on $S$, such that
for every complex Hilbert space $\Hh$ and every
hermitian $F\in\B(\Hh)[x]$ with $\deg F\le2$,
\[
F\succeq0\text{ on }K
\quad\Longrightarrow\quad
F\in\M_{\Hh}(S)_{\le4\ell-2+\deg\sigma_1+\eta_S}.
\]
\end{corollary}

\begin{proof}
Since $\Delta\ge0$ on $K$ and $M(S)$ is saturated,
there exists a representation
$\Delta=\sum_{j=0}^s g_j\tau_j$, where $g_0=1$
and each $\tau_j\in\SOS_{\C}$.
Let $\vd_S$ denote the weighted degree of this certificate
and set $\eta_S=\vd_S-2 \ell$. Since $\deg\Delta=2\ell$, we have $\eta_S\geq 0.$

Let $\Hh$ be a complex Hilbert space and let
$F\in\B(\Hh)[x]$ be hermitian, with $\deg F\le2$
and $F\succeq0$ on $K$.
By Theorem~\ref{thm:intervals}, there exist
$\Sigma_0,\Sigma_1\in\SOS_{\Hh}$ such that
$F=\Sigma_0+\Delta\Sigma_1$, with degree bounds
$\deg\Sigma_0,\deg(\Delta\Sigma_1)\le4\ell-2+\deg\sigma_1$.
Substituting the scalar certificate for $\Delta$ gives
\[
F=\Sigma_0+\sum_{j=0}^s g_j\tau_j\Sigma_1.
\]
Since each $\tau_j\Sigma_1\in\SOS_{\Hh}$,
the resulting certificate has weighted degree at most
$4\ell-2+\deg\sigma_1+\eta_S$. Thus,
$
F\in\M_{\Hh}(S)_{\le4\ell-2+\deg\sigma_1+\eta_S},
$
as claimed.
\end{proof}

\subsection{Explicit multiplier}
\label{sssec:explicit-formula}
Recall that
\[
P(x)=\prod_{j=1}^{\ell}(x-a_j),
\qquad
Q(x)=\prod_{j=1}^{\ell}(x-b_j).
\]
By \eqref{eq:v^Tu} and the diagonal form of $D$,
\[
P(x)-Q(x)
=
\sum_{j=1}^{\ell}v_j^2
\prod_{k\ne j}(x-b_k).
\]
Evaluating at $x=b_j$ gives
\[
v_j^2=\frac{P(b_j)}{Q'(b_j)}>0,
\qquad j=1,\ldots,\ell.
\]
Thus, using \eqref{eq:u} the polynomial $w=u^*u$ is given by
\[
w(x)=
\sum_{j=1}^{\ell}
\frac{P(b_j)}{Q'(b_j)}
\prod_{k\ne j}(x-b_k)^2.
\]

The following lemma will be used in Section~\ref{sec:uniform} to establish uniform degree bounds for operator positivity certificates over families of compact interval unions contained in a fixed compact interval, with component lengths and separating gaps bounded below by positive constants.

\begin{lemma}\label{lem:h-continuity}
Let $K$ and $\ell$ be as in \eqref{eq:endpoint-data},
with $\ell\ge1$ and pairwise distinct endpoints
$a_1,b_1,\ldots,a_\ell,b_\ell$.
Then the coefficients of $w$ depend continuously
on these endpoints.
\end{lemma}

\begin{proof}
The coefficients of $w$ are rational functions
of the endpoints, with denominators
\[
Q'(b_j)=\prod_{k\ne j}(b_j-b_k),
\]
which are nonzero since the $b_j$ are distinct.
\end{proof}

\section{Proof of Theorem~\ref{thmA}}
\label{sec:points}

We extend the positivity certificates of Theorem~\ref{thm:intervals}
from finite unions of nondegenerate compact intervals to compact semialgebraic sets that also contain isolated points. 

Let
\begin{equation}\label{def:J}
 K_0=\bigcup_{j=1}^{\ell}[a_j,b_j],\qquad
 a_1<b_1<\cdots<a_\ell<b_\ell,\qquad \ell\ge1,
\end{equation}
and recall the defining polynomial
\begin{equation}\label{def:Delta}
 \Delta(x)=-P(x)Q(x)
 =-\prod_{j=1}^{\ell}(x-a_j)(x-b_j),
 \qquad K_0=\{x\in\R:\Delta(x)\ge0\}.
\end{equation}

\subsection{Adjoining a single isolated point}
Fix $\xi\in \R\setminus K_0$ and put $\delta=-\Delta(\xi)>0$.
Define
\begin{equation}\label{eq:point-matrix}
 G_\xi(x)=
 \begin{pmatrix}
 \delta+\Delta(x)&(x-\xi)\Delta(x)\\
 (x-\xi)\Delta(x)&(x-\xi)^2\Delta(x)
 \end{pmatrix}.
\end{equation}

The following lemma describes the positivity set of $G_\xi$.

\begin{lemma}\label{lem:point-matrix}
For $\xi\in\R\setminus K_0$,
let $G_\xi$ be the $2\times2$ matrix polynomial defined by
\eqref{eq:point-matrix}. Then $G_\xi(\xi)=0$ and
\[
 \{x\in\R \ \colon \ G_\xi(x)\succeq0\} \,=\, K_0\cup\{\xi\}.
\]
\end{lemma}
\begin{proof}
The identity
\[
 G_{\xi}(x)=\delta\begin{pmatrix}1&0\\0&0\end{pmatrix}
 +\Delta(x)\begin{pmatrix}1\\x-{\xi}\end{pmatrix}
           \begin{pmatrix}1&x-{\xi}\end{pmatrix}
\]
gives $G_{\xi}(x)\succeq0$ for $x\in K_0$.
At $x={\xi}$, the equality $\Delta({\xi})=-\delta$ yields $G_{\xi}({\xi})=0$.
Thus 
\[
K_0\cup\{{\xi}\}\subseteq\{x\in\R:G_{\xi}(x)\succeq0\}.
\]

Conversely, if $G_{\xi}(x)\succeq0$, then
\[
\det G_{\xi}(x)=\delta\Delta(x)(x-\xi)^2\ge0.
\]
If $x\ne\xi$, then $\delta(x-\xi)^2>0$, so
$\Delta(x)\ge0$ and hence $x\in K_0$.
\end{proof}

Positivity at $\xi$ allows a factorization through $G_\xi$ without increasing the degrees of the operator columns.
\begin{lemma}\label{lem:point-matching}
Fix $\xi\in\R\setminus K_0$, and let $G_\xi$
be defined by \eqref{eq:point-matrix}.
Suppose that $F\in\B(\Hh)[x]$ admits a representation
\[
F=A^*A+\Delta B^*B,
\]
where $A,B\in\B(\Hh)[x]$ satisfy
$\deg A,\deg B\le\vd$ for some integer $\vd\ge0$.
\smallskip

If $F(\xi)\succeq0$, then there exist operator polynomials
$L_\xi,W_\xi,E_\xi\in\B(\Hh)[x]$ such that, defining
$V_\xi=\col(W_\xi,E_\xi)$, we have
$\deg L_\xi,\deg V_\xi\le\vd$ and
\begin{equation}\label{eq:point-factorization}
F=L_\xi^*L_\xi+
V_\xi^*(I_{\Hh}\otimes G_\xi)V_\xi.
\end{equation}
\end{lemma}

\begin{proof}
Since $\Delta({\xi})=-\delta$, the assumption $F({\xi})\succeq0$ gives
\[
 \delta B({\xi})^*B({\xi})\preceq A({\xi})^*A({\xi}).
\]
By Douglas' factorization lemma~\cite{Doug66},
there exists a contraction $\Gamma:\Hh\to\Hh$ such that
\[
 \Gamma A({\xi})=\sqrt \delta\,B({\xi}).
\]
Define
\[
 \begin{aligned}
 W_{\xi}(x)&=\frac{\Gamma A(x)}{\sqrt \delta},
 & \qquad E_{\xi}(x)&=\frac{B(x)- W_{\xi}(x)}{x-{\xi}},\\
 L_{\xi}(x)&=(I_{\Hh_A}-\Gamma^*\Gamma)^{1/2}A(x),
 & \qquad V_{\xi}(x)&=\col(W_{\xi}(x),E_{\xi}(x)).
 \end{aligned}
\]
Because $W_{\xi}({\xi})=B({\xi})$, the column $E_{\xi}$ is a polynomial,
with $\deg E_{\xi}\le \vd-1$.
Consequently, $\deg L_{\xi},\deg V_{\xi}\le \vd$.

\medskip

The definitions of $L_\xi$ and $W_\xi$ give
\[
L_\xi^*L_\xi
=A^*(I_{\Hh_A}-\Gamma^*\Gamma)A
=A^*A-\delta W_\xi^*W_\xi.
\]
Moreover, $B=W_\xi+(x-\xi)E_\xi$. Hence, by
\eqref{eq:point-matrix},
\begin{align*}
V_\xi^*(I_{\Hh_B}\otimes G_\xi)V_\xi
&=\delta W_\xi^*W_\xi
+\Delta\bigl(W_\xi+(x-\xi)E_\xi\bigr)^*
       \bigl(W_\xi+(x-\xi)E_\xi\bigr)\\
&=\delta W_\xi^*W_\xi+\Delta B^*B.
\end{align*}
Combining these identities gives
\[
L_\xi^*L_\xi+
V_\xi^*(I_{\Hh_B}\otimes G_\xi)V_\xi
=A^*A+\Delta B^*B=F,
\]
which proves \eqref{eq:point-factorization}.
\end{proof}

\subsection{Adjoining finitely many isolated points}
Let $\xi_1,\ldots,\xi_r\in\R\setminus K_0$ be distinct,
with $r\ge1$, and put $K=K_0\cup\{\xi_1,\ldots,\xi_r\}$.
The following lemma constructs interpolation polynomials
for the points $\xi_1,\ldots,\xi_r$ whose squared sum
is bounded by $1$ on $K_0$.

\begin{lemma}\label{lem:point-cutoffs}
There are real polynomials $\psi_1,\ldots,\psi_r$, depending only on
$K$, such that
\[
 \psi_i(\xi_j) = \delta_{ij},\qquad
 \sum_{i=1}^r \psi_i(x)^2\le1\quad(x\in K_0).
\]
\end{lemma}

\begin{proof}
Let $\ell_i$ be the Lagrange interpolation polynomials
defined in \eqref{eq:Lagrange}. For each $i$, choose
$n_i>0$ such that $|x-\xi_i|<n_i$ for all $x\in K_0$.
For $N\ge1$, set
\[
\psi_{i,N}(x)=\ell_i(x)
\left(1-\frac{(x-\xi_i)^2}{n_i^2}\right)^N.
\]
The interpolation conditions hold for every $N$.
Since $\psi_{i,N}\to0$ uniformly on $K_0$ as
$N\to\infty$, choosing $N$ sufficiently large
gives the desired inequality.
\end{proof}

The interpolation polynomials allow us
to adjoin all isolated points with a fixed additive increase in the
weighted certificate degree.

\begin{proposition}\label{prop:points}
Let $\xi_1,\ldots,\xi_r$ be
distinct points in $\R\setminus K_0$, with $r\ge1$, and set
$K=K_0\cup\{\xi_1,\ldots,\xi_r\}$.
Let $S\subset\R[x]$ be finite with $K_S=K$ and $M(S)$ saturated.
Then there exists a nonnegative integer $\eta_S$, depending only on $S$, such that, for every 
Hilbert space $\Hh$, every nonnegative integer $n$, and every hermitian
$F\in\B(\Hh)[x]$,
\[
 F\in\M_{\Hh}(\{\Delta\})_{\le n},\quad
 F(\xi_i)\succeq0\ (i=1,\ldots,r)
 \quad\Longrightarrow\quad
 F\in\M_{\Hh}(S)_{\le n+\eta_S}.
\]
\end{proposition}

\begin{proof}
Let $\psi_1,\ldots,\psi_r$ be the interpolation polynomials
from Lemma~\ref{lem:point-cutoffs}, and set
\[
a=1-\sum_{i=1}^r\psi_i^2,\qquad
\Psi_i=\psi_i^2G_{\xi_i}\quad(1\le i\le r),
\]
where $G_{\xi_i}$ is defined in \eqref{eq:point-matrix}.
Since $a\ge0$ on $K_0$ and $a(\xi_j)=0$ for every $j$,
both $a$ and $a\Delta$ are nonnegative on $K$.
By Lemma~\ref{lem:point-matrix}, each $\Psi_i$ is
positive semidefinite on $K_0$.
Moreover, the interpolation conditions imply
$\Psi_i(\xi_j)=0$ for all $i,j$.
Thus $\Psi_i\succeq0$ on $K$.
\medskip

Since $M(S)$ is saturated, both $a$ and $a\Delta$
admit scalar certificates in $M(S)$.
By \cite[Theorem~2.1]{Zalar}, each $\Psi_i$
admits a matrix certificate in $\M_{\C^2}(S)$.
Choose a nonnegative integer $\eta_S$ large enough
to bound the weighted degrees of all these certificates
and such that $2\deg\psi_i\le\eta_S$ for every $i$.
\medskip

Let $\Hh$ be a complex Hilbert space and let
$F\in\B(\Hh)[x]$ satisfy the assumptions of the proposition.
Since $F\in\M_{\Hh}(\{\Delta\})_{\le n}$,
there exist $A,B\in\B(\Hh)[x]$ such that
$F=A^*A+\Delta B^*B$ has weighted degree at most $n$.
Set $\widetilde n=\lfloor n/2\rfloor$.
Then $\deg A,\deg B\le\widetilde n$.
Applying Lemma~\ref{lem:point-matching} at each $\xi_i$
gives polynomials $L_i,W_i,E_i\in\B(\Hh)[x]$,
with $V_i=\col(W_i,E_i)$ and
$\deg L_i,\deg V_i\le\widetilde n$, such that
\[
F=L_i^*L_i+
V_i^*(I_{\Hh}\otimes G_{\xi_i})V_i.
\]
Using $1=a+\sum_{i=1}^r\psi_i^2$, we obtain
\begin{align}
F
&=aF+\sum_{i=1}^r\psi_i^2F \nonumber\\
&=aA^*A+(a\Delta)B^*B
  +\sum_{i=1}^r(\psi_iL_i)^*(\psi_iL_i)
  +\sum_{i=1}^r
  V_i^*(I_{\Hh}\otimes\Psi_i)V_i.
\label{eq:all-points}
\end{align}

Substituting the chosen scalar certificates for $a$
and $a\Delta$, together with the matrix certificates
for $\Psi_i$, into \eqref{eq:all-points}, we obtain a
certificate for $F$ in $\M_{\Hh}(S)$.
Here \ref{item:global-factorization} of
Proposition~\ref{prop:factorization} allows us to
express the resulting positive semidefinite operator
polynomial multipliers as hermitian squares without
increasing their degrees.
Since $\deg A,\deg B,\deg L_i,\deg V_i\le\widetilde n$,
the resulting certificate has weighted degree at most
$\eta_S+2\widetilde n\le n+\eta_S$.
Thus
$
F\in\M_{\Hh}(S)_{\le n+\eta_S}.
$
\end{proof}

\subsection{Completion of proof of Theorem~\ref{thmA}}
\begin{proof}[Proof of Theorem~\ref{thmA}]
We first work over a complex Hilbert space $\Hh$. If $K$ is finite, Proposition~\ref{prop:finite-set} gives the required bound with $c_S=\eta_S$. We therefore assume that $K$ contains at least one nondegenerate interval.

\smallskip

Let $K_0$ be the union of all nondegenerate interval
components of $K$, and let $\ell$ be their number.
Apply Theorem~\ref{thm:intervals} to $K_0$, fix the
corresponding polynomial $\sigma_1$, and put
$\beta=4\ell-2+\deg\sigma_1$.
Let $\mathcal E$ be an arbitrary complex Hilbert space.
By Theorem~\ref{thm:intervals},
Corollary~\ref{cor:interval-transfer}, and
Proposition~\ref{prop:points}, there exists a
nonnegative integer $\eta_S$, depending only on $S$,
such that every hermitian $F\in\B(\mathcal E)[x]$
with $\deg F\le2$ and $F\succeq0$ on $K$ satisfies
\[
F\in\M_{\mathcal E}(S)_{\le N_S},
\qquad N_S=\beta+\eta_S.
\]
In particular, $N_S$ is independent of $\mathcal E$.

\smallskip

By Theorem~\ref{thm:quadratic-reduction}, for every $\vd\ge1$, it follows that
\[
F\succeq0\text{ on }K,\qquad \deg F\le\vd
\quad\Longrightarrow\quad
F\in\M_{\Hh}(S)_{\le N_S+2\lceil\vd/2\rceil-2}
\subseteq\M_{\Hh}(S)_{\le\vd+N_S-1}.
\]
Hence the required degree bound holds with $c_S=N_S-1$.  For $\vd=0$, the positive constant operator $F$ is a single hermitian square. The same bounds for real Hilbert spaces follow from Lemma~\ref{lem:real-descent}. The converse follows by evaluating the certificates on $K$.
\end{proof}

\subsection{Computation of the additive constant}
\label{subsec:degree-constants}

We now discuss the additive constants obtained in the proof of Theorem~\ref{thmA}. Suppose $K$ contains $\ell\ge1$
nondegenerate intervals. Choose $\sigma_0,\sigma_1\in\SOS_{\C}$ satisfying
\eqref{eq:scalar2}.
Let $\eta_S$ be the additive degree increase from
Corollary~\ref{cor:interval-transfer} or
Proposition~\ref{prop:points}, as appropriate.
Then the degree bound holds with
\[
c_S=4\ell-3+\deg\sigma_1+\eta_S.
\]
If $K$ is finite, Proposition~\ref{prop:finite-set}
gives $c_S=\eta_S$.
These constants need not be optimal.

\medskip

Although $\Delta$ and $w$ are explicit,
the degrees of the required SOS multipliers are
not known in general. If the generators in $S$
have real algebraic coefficients, admissible constants
can be computed by increasing the degree bounds
and testing semidefinite feasibility. These feasibility conditions are
decidable by quantifier elimination over the real
algebraic numbers, and saturation guarantees
termination. The interpolation polynomials for
isolated points can also be constructed effectively.

\begin{remark}\label{rem:explicit-computations}
Explicit computations illustrating the preceding construction,
together with Mathematica verification code, are available in
the accompanying repository~\cite{Code}. For the natural
descriptions $S_0$ and $S_1$ of
$K_0=[0,1]\cup[2,3]$,
$K_1=K_0\cup\{4\},$
respectively, the construction gives admissible additive constants
$c_{S_0}=12$ and $c_{S_1}=18$ in Theorem~\ref{thmA}.
These constants are not asserted to be optimal.
\end{remark}

\section{Dependence on the underlying set}\label{sec:nonuniform}
The additive constant in Theorem~\ref{thm:main} depends on the compact set and its scalar description. The following examples show that the additive constants can become arbitrarily large, even for quadratic $2\times2$ matrix polynomials on unions of two compact intervals. In the first example, the right endpoint of one interval tends to infinity. In the second, both intervals remain in $[0,1]$, but one of them shrinks to a point.

\medskip

Both examples use the following compactness property of matrix
polynomials with uniformly bounded degree.

\begin{lemma}\label{lem:compactness}
Let $I\subset\R$ be a nondegenerate compact interval,
and fix integers $m\ge1$ and $\vd\ge0$.
Suppose $(F_k)$ is a sequence in $\C[x]^{m\times m}$
with $\deg F_k\le\vd$ for every $k$, satisfying
\begin{equation}\label{eq:compactness-bound}
\sup_k\sup_{x\in I}\|F_k(x)\|_2<\infty,
\end{equation}
where $\|\cdot\|_2$ denotes the spectral norm.
Then $(F_k)$ has a subsequence that converges
coefficientwise to a polynomial $F\in\C[x]^{m\times m}$
of degree at most $\vd$.
Moreover, if each $F_k$ is positive semidefinite
on $\R$, then so is $F$.
\end{lemma}

\begin{proof}
The polynomials in $\C[x]^{m\times m}$ of degree at most
$\vd$ form a finite-dimensional vector space.
For $G(x)=\sum_{j=0}^{\vd}A_jx^j$, define
\[
\|G\|_{\mathrm{coeff}}
:=\max_{0\le j\le\vd}\|A_j\|_2,
\qquad
\|G\|_I:=\sup_{x\in I}\|G(x)\|_2.
\]
Since $I$ is nondegenerate, a polynomial vanishing on $I$
must vanish identically. Applied entrywise, this shows that
$\|\cdot\|_I$ is a norm. The two norms are therefore equivalent
by finite dimensionality.
Hence \eqref{eq:compactness-bound} yields a coefficientwise
convergent subsequence, which converges at every
$x\in\R$. Since the cone
of positive semidefinite matrices is closed, the limit is
positive semidefinite on $\R$.
\end{proof}

\subsection{Example I} 
Consider the quadratic matrix polynomial
\begin{equation}\label{eq:thesis-example}
F(x) \,=\, 
\begin{pmatrix}
x+1 \ &\ \sqrt6\\
\sqrt6 \ &\ x^2-4x+6
\end{pmatrix}
\end{equation} 
and the scalar polynomial $q(x)=(x-1)(x-2)$. The following nonmembership result is a special case of \cite[Proposition~3.4, Claim~1]{Zalar}.

\begin{lemma}\label{lem:noncompact-obstruction}
Let $F$ be as in \eqref{eq:thesis-example} and
$q=(x-1)(x-2)$. Then
\[
F\succeq0\text{ on }[0,1]\cup[2,\infty),
\qquad
F\notin\mathcal T_{\C^2}(\{x,q\}).
\]
\end{lemma}

\begin{proof}
The diagonal entries of $F$ are positive on
$[0,1]\cup[2,\infty)$, and
$\det F(x)=x(x-1)(x-2)\ge0$ on this set.
Hence $F\succeq0$ on $[0,1]\cup[2,\infty)$.
The assertion $F\notin\mathcal T_{\C^2}(\{x,q\})$
follows from \cite[Proposition~3.4, Claim~1]{Zalar}.
\end{proof}

Restricting $F$ to compact truncations gives certificates whose
minimum weighted degree must diverge as the truncation expands.

\begin{proposition}\label{prop:truncations}
For $n\ge3$, let $K_n=[0,1]\cup[2,n]$ and
$S_n=\{x,q,1-x/n\}$. Then, for $F$ as in
\eqref{eq:thesis-example},
\[
\min\bigl\{k\ge0:
F\in\mathcal T_{\C^2}(S_n)_{\le k}\bigr\}
\longrightarrow\infty
\qquad(n\to\infty).
\]
\end{proposition}

\begin{proof}
For each $n\ge3$, Theorem~\ref{thm:scalar} shows that
$S_n$ is a saturated scalar description of $K_n$.
Since $F\succeq0$ on $K_n$, Theorem~\ref{thm:main}
implies that
\[
F\in\M_{\C^2}(S_n)\subseteq\mathcal T_{\C^2}(S_n).
\]

Suppose, for a contradiction, that there exist an integer
$\vd\ge0$ and a sequence $n_j\to\infty$ such that
$F\in\mathcal T_{\C^2}(S_{n_j})_{\le\vd}$ for every $j$.
Then
\[
F=\sum_{a,b,c\in\{0,1\}}
x^a q^b(1-x/n_j)^c\Sigma_{j,abc},
\qquad \Sigma_{j,abc}\in\SOS_{\C^2},
\]
where each representation has weighted degree at most $\vd$.
On $I=[1/4,3/4]$, all eight scalar weights are bounded below
by a constant $\delta>0$ independent of $j$.
Since every summand is positive semidefinite on $I$, we obtain
\[
0\preceq\Sigma_{j,abc}(x)\preceq\delta^{-1}F(x)
\qquad(x\in I).
\]
Thus the polynomials $\Sigma_{j,abc}$ are uniformly bounded
in spectral norm on $I$, and $\deg\Sigma_{j,abc}\le\vd$.
\medskip

Applying Lemma~\ref{lem:compactness} successively to the eight
sequences, we obtain a common subsequence along which
$\Sigma_{j,abc}$ converges coefficientwise to $\Sigma_{abc}$.
Each limit is positive semidefinite on $\R$ and therefore belongs
to $\SOS_{\C^2}$ by \ref{item:global-factorization} of
Proposition~\ref{prop:factorization}.
Since $1-x/n_j$ converges coefficientwise to $1$, passing to
the limit yields
\[
 F=
 \sum_{a,b\in\{0,1\}}
 x^a q^b(\Sigma_{ab0}+\Sigma_{ab1}).
\]
Thus $F\in\mathcal T_{\C^2}(\{x,q\})$, contradicting Lemma~\ref{lem:noncompact-obstruction}.
\end{proof}
\subsection{Example II} \label{ssec:exampleII}
We next construct a family of unions of two disjoint compact intervals for which the gap
stays bounded away from zero while one interval component
collapses. 

\medskip

For $0<\varepsilon<1/2$, define
\[
 \begin{gathered}
 K_\varepsilon=[0,\varepsilon^2]\cup[1/2,1],
 \qquad
 q_\varepsilon (x) =(x-\varepsilon^2)(x-1/2),\\
 \theta(x) = x(1-x),\qquad
 \rho_\varepsilon=1+2\varepsilon-\varepsilon^3,
 \end{gathered}
\]
and let
\begin{equation}\label{eq:counterexample}
 F_\varepsilon (x) =
 \begin{pmatrix}
 \rho_\varepsilon(x^2+\varepsilon \theta (x))
   &\rho_\varepsilon \theta (x) \\
 \rho_\varepsilon \theta (x)
   &\varepsilon(1-x)^2+(1-\varepsilon^2) \theta (x)
 \end{pmatrix}.
\end{equation}

\begin{lemma}\label{lem:family-positive}
The quadratic matrix polynomial $F_\varepsilon$ is positive
semidefinite on $K_\varepsilon$ for every
$0<\varepsilon<1/2$.
\end{lemma}

\begin{proof}
Since $\rho_\varepsilon>0$, both diagonal entries are nonnegative
on $[0,1]$.  The determinant of $F_{\varepsilon}$ is given by
$ \det F_\varepsilon
 =2\rho_\varepsilon \theta q_\varepsilon. $
The right-hand side is nonnegative on $K_\varepsilon$, proving
the assertion.
\end{proof}

We include a proof of the following elementary fact.
\begin{lemma}\label{lem:double-zero}
Let $H\in\C[x]^{m\times m}$ be hermitian and positive semidefinite
on $\R$. If $H(0)=0$, then $H'(0)=0$, and $x^2$ divides $H$
entrywise.
\end{lemma}

\begin{proof}
For every $\xi\in\C^m$, the scalar polynomial
$\xi^*H(x)\xi$ is nonnegative on $\R$ and vanishes at zero.
Its derivative therefore vanishes there, so
$
 \xi^*H'(0)\xi=0$.
Polarization gives $H'(0)=0$. Together with $H(0)=0$, this
implies the asserted divisibility.
\end{proof}

\subsection{Proof of Theorem~\ref{thmB}}

\begin{proof}[Proof of Theorem~\ref{thmB}]
For $0<\varepsilon<1/2$, define
$K_\varepsilon=[0,\varepsilon^2]\cup[1/2,1]$
and $q_\varepsilon=(x-\varepsilon^2)(x-1/2)$,
and set $S_{K_\varepsilon}=\{x,1-x,q_\varepsilon\}$.
By Theorem~\ref{thm:scalar}, $S_{K_\varepsilon}$ is a
saturated scalar description of $K_\varepsilon$.
For the matrix polynomial $F_\varepsilon$ defined in
\eqref{eq:counterexample}, Lemma~\ref{lem:family-positive}
gives $F_\varepsilon\succeq0$ on $K_\varepsilon$.
Hence Theorem~\ref{thm:main} yields
\[
F_\varepsilon\in\M_{\C^2}(S_{K_\varepsilon})
\subseteq\mathcal T_{\C^2}(S_{K_\varepsilon}).
\]

Suppose, for a contradiction, that there exist an integer
$\vd\ge0$ and a sequence $\varepsilon_j\downarrow0$ such that
$F_{\varepsilon_j}\in\mathcal T_{\C^2}(S_{K_{\varepsilon_j}})_{\le\vd}$
for every $j$. Then
\begin{equation}\label{eq:bounded-preordering}
F_{\varepsilon_j}(x)
=\sum_{a,b,c\in\{0,1\}}
x^a(1-x)^bq_{\varepsilon_j}^c\Sigma_{j,abc},
\qquad \Sigma_{j,abc}\in\SOS_{\C^2},
\end{equation}
where each representation has weighted degree at most $\vd$.
On $I=[3/4,7/8]$, all eight scalar weights are bounded below
by a constant $\delta>0$ independent of $j$.
Since every summand is positive semidefinite on $I$, we obtain
\[
0\preceq\Sigma_{j,abc}(x)
\preceq\delta^{-1}F_{\varepsilon_j}(x)
\qquad(x\in I).
\]
Since the polynomials $F_{\varepsilon_j}$ are uniformly bounded
in spectral norm on $I$, the same holds for all $\Sigma_{j,abc}$.
Moreover, $\deg\Sigma_{j,abc}\le\vd$.

\medskip

Applying Lemma~\ref{lem:compactness} successively to the
eight sequences $(\Sigma_{j,abc})_j$, indexed by
$(a,b,c)\in\{0,1\}^3$, we obtain a common subsequence
along which $\Sigma_{j,abc}$ converges coefficientwise
to $\Sigma_{abc}$ for every $(a,b,c)\in\{0,1\}^3$.
Each limit $\Sigma_{abc}$ is positive semidefinite on $\R$.

\medskip

We first extract additional information by evaluating the
certificates at zero. We have
\[
 F_{\varepsilon_j}(0)=
 \begin{pmatrix}0&0\\0&\varepsilon_j\end{pmatrix},
 \qquad
 q_{\varepsilon_j}(0)=\varepsilon_j^2/2>0.
\]
For $e_1=(1,0)^{\mathsf T}$, equation
\eqref{eq:bounded-preordering} gives
\[ 0=
 \sum_{b,c\in\{0,1\}}
 q_{\varepsilon_j}(0)^c
 e_1^{\mathsf T}\Sigma_{j,0bc}(0)e_1.
\]
Every summand is nonnegative and every scalar weight is positive,
so $e_1^{\mathsf T}\Sigma_{j,0bc}(0)e_1=0$ for all $b,c$.
Positive semidefiniteness then implies
$\Sigma_{j,0bc}(0)e_1=0$. Passing to the limit yields
\begin{equation}\label{eq:kernel-limit}
 \Sigma_{0bc}(0)e_1=0
 \qquad(b,c\in\{0,1\}).
\end{equation}

The coefficientwise limits of $F_{\varepsilon_j}$ and
$q_{\varepsilon_j}$ are
\[
 F_0(x)=xG_0(x),\qquad
 G_0(x)=
 \begin{pmatrix}
 x&1-x\\
 1-x&1-x
 \end{pmatrix},
 \qquad q_0(x)=x(x-1/2).
\]
Passing to the limit in \eqref{eq:bounded-preordering} gives
\begin{equation}\label{eq:lim-rep}
 F_0(x)=
 \sum_{a,b,c\in\{0,1\}}
 x^a(1-x)^bq_0(x)^c\Sigma_{abc}.
\end{equation}
Evaluating at zero, we obtain
$ \Sigma_{000}(0)+\Sigma_{010}(0)=0. $
Both matrices are positive semidefinite, so each is zero.
Lemma~\ref{lem:double-zero} therefore shows that $x^2$ divides
both $\Sigma_{000}$ and $\Sigma_{010}$ entrywise.

\medskip

Thus every term in the limiting representation \eqref{eq:lim-rep} is
divisible by $x$. After dividing by $x$ and setting
$x=0$, the terms with $a=c=0$ vanish by the
divisibility established above, while those with $a=c=1$ retain a factor of $x$. Hence
\[
 \begin{pmatrix}0&1\\1&1\end{pmatrix}
 =
 \Sigma_{100}(0)+\Sigma_{110}(0)
 -\frac12\bigl(\Sigma_{001}(0)+\Sigma_{011}(0)\bigr).
\]
Set $A=\Sigma_{100}(0)+\Sigma_{110}(0)\succeq0$.
Since $\Sigma_{001}(0)e_1=\Sigma_{011}(0)e_1=0$
by \eqref{eq:kernel-limit}, the preceding identity gives
$ Ae_1= (0,1)^{\mathsf T}$.
But $e_1^{\mathsf T} A e_1=0$ and $A\succeq0$ imply $Ae_1=0$, a
contradiction. Thus the minimum weighted certificate degree
tends to infinity as $\varepsilon\downarrow0$.
\end{proof}

\section{Proof of Theorem~\ref{thmC}}
\label{sec:uniform}

We prove Theorem~\ref{thmC} in this section. 
Fix a compact interval $J\subset\R$ and $\mu,\omega>0$, as in Theorem~\ref{thmC}. Since there is nothing to prove if $\mathcal K(J,\mu,\omega)$ is empty, we assume that this family is nonempty. Let $|J|$ denote the length of $J$. Since every $K\in\mathcal K(J,\mu,\omega)$ with $\ell$ components satisfies
\[
\ell\mu+(\ell-1)\omega\le |J|,
\]
only finitely many values of $\ell$ occur. For each such $\ell$, let $\mathcal P_\ell\subset J^{2\ell}$ be the compact set of endpoint tuples
$\mathbf e=(a_1,b_1,\ldots,a_\ell,b_\ell)$ corresponding to sets in $\mathcal K(J,\mu,\omega)$ with $\ell$ components.

\medskip

For each $\mathbf e=(a_1,b_1,\ldots,a_\ell,b_\ell)\in\mathcal P_\ell$, let $K_{\mathbf e}$ denote the corresponding union of intervals and, recalling \eqref{eq:P-Q-r-Delta}, set $\Delta_{\mathbf e}(x)=-\prod_{j=1}^{\ell}(x-a_j)(x-b_j)$. 
Let $w_{\mathbf e}$ be obtained from
the rank-one pencil associated with
$K_{\mathbf e}$.

\medskip

We first show that the degrees of the SOS
multipliers in the identity
\begin{equation}\label{eq:uniform-scalar}
1=w_{\mathbf e}\sigma+
\Delta_{\mathbf e}\tau,
\qquad \sigma,\tau\in\SOS_{\C},
\end{equation}
can be bounded independently of
$\mathbf e\in\mathcal P_\ell$.

\begin{lemma}\label{lem:strict positive}
For every $\mathbf e\in\mathcal P_\ell$,
there exist $\sigma,\tau\in\SOS_{\C}$
satisfying \eqref{eq:uniform-scalar}
and a number $\varepsilon>0$ such that
\[
\sigma(x),\tau(x)\ge\varepsilon
\qquad(x\in\R).
\]
\end{lemma}

\begin{proof}
Since $K_{\mathbf e}$ is compact,
we can choose $\varepsilon>0$ such that
$ f=1-\varepsilon
(w_{\mathbf e}+\Delta_{\mathbf e}) $
is nonnegative on $K_{\mathbf e}$.
Since $w_{\mathbf e}>0$ on $\R$ and $\Delta_{\mathbf e}$ has simple zeros at all the endpoints of $K_{\mathbf e}$, the quadratic module $M(\{w_{\mathbf e}\Delta_{\mathbf e}\})$ is saturated by Theorem~\ref{thm:scalar}.
Hence, as
$\Delta_{\mathbf e}f\ge0$ on $K_{\mathbf e}$,
there exist
$\widetilde\sigma_0,\widetilde\sigma_1
\in\SOS_{\C}$ such that
\[
\Delta_{\mathbf e}f
=
\widetilde\sigma_0+
w_{\mathbf e}\Delta_{\mathbf e}
\widetilde\sigma_1.
\]
At every endpoint of $K_{\mathbf e}$,
the polynomial $\widetilde\sigma_0$
vanishes. Since it is nonnegative on
$\R$, each such zero has even
multiplicity. Therefore,
$ \widetilde\sigma_0 =\Delta_{\mathbf e}^2\widetilde\tau $
for some $\widetilde\tau\in\SOS_{\C}$,
since the quotient is nonnegative on $\R$.
Dividing by $\Delta_{\mathbf e}$ gives
\[
f=\Delta_{\mathbf e}\widetilde\tau + w_{\mathbf e}\widetilde\sigma_1.
\]
Thus, 
$ 1=w_{\mathbf e}
(\widetilde\sigma_1+\varepsilon)
+\Delta_{\mathbf e}
(\widetilde\tau+\varepsilon). $
Taking
$ \sigma=\widetilde\sigma_1+\varepsilon $
and
$\tau=\widetilde\tau+\varepsilon $
proves the lemma.
\end{proof}

For $N\in\mathbb N$, let
\[
V_N=
\R[x]_{\le2\ell+2N}
\oplus
\R[x]_{\le2\ell-2+2N},
\qquad
W_N=\R[x]_{\le4\ell-2+2N}.
\]
For $\mathbf e\in\mathcal P_\ell$, define
the linear map
\[
\mathcal L_{\mathbf e,N}:V_N
\longrightarrow W_N; \qquad 
(f,g) \longmapsto
w_{\mathbf e}f+\Delta_{\mathbf e}g.
\]

\begin{lemma}\label{lem:L-surjective}
For every $\mathbf e\in\mathcal P_\ell$
and every $N\in\mathbb N$, the linear map
$\mathcal L_{\mathbf e,N}$ is surjective.
\end{lemma}

\begin{proof}
Since $w_{\mathbf e}>0$ on $\R$ and
all zeros of $\Delta_{\mathbf e}$ are real,
the polynomials $w_{\mathbf e}$ and
$\Delta_{\mathbf e}$ are coprime.
By B\'ezout's identity, there exist
$a,b\in\R[x]$ such that
\[aw_{\mathbf e}+b\Delta_{\mathbf e}=1.\]

Let $q\in W_N$. Multiplying this identity by $q$
and dividing $qa$ by $\Delta_{\mathbf e}$, we obtain
$qa=\widetilde f\Delta_{\mathbf e}+f$,
where $\deg f\le2\ell-1$. Thus,
\begin{equation}\label{eq:surjective-representation}
q=w_{\mathbf e}f+\Delta_{\mathbf e}g,
\qquad g=w_{\mathbf e}\widetilde f+qb.
\end{equation}
Since $\deg w_{\mathbf e}=2\ell-2$, we have
$\deg(w_{\mathbf e}f)\le4\ell-3$.
Together with $\deg q\le4\ell-2+2N$ and
$\deg\Delta_{\mathbf e}=2\ell$, this implies
$\deg g\le2\ell-2+2N$.
Thus $(f,g)\in V_N$, proving the surjectivity
of $\mathcal L_{\mathbf e,N}$.
\end{proof}

\begin{remark}
The polynomial $f$ in \eqref{eq:surjective-representation}
satisfies the stronger degree bound $\deg f\le2\ell-1$,
whereas the first component of $V_N$ allows degree
up to $2\ell+2N$. This additional degree allowance
is not needed for surjectivity. It is used in
Lemma~\ref{lem:strict positive plus degrees},
where the first SOS multiplier is perturbed by
$-\delta\Delta_{\mathbf e}(1+x^2)^N$, a polynomial
of degree $2\ell+2N$. Unlike the polynomial $f$
in \eqref{eq:surjective-representation}, the
resulting multiplier must be strictly positive
and have the prescribed degree, ensuring that
it lies in the interior of the corresponding
truncated SOS cone.
\end{remark}

\begin{lemma}\label{lem:strict positive plus degrees}
Fix $\mathbf e\in\mathcal P_\ell$,
and let $\sigma,\tau\in\SOS_{\C}$
be as in Lemma~\ref{lem:strict positive}.
Suppose $N\in\mathbb N$ satisfies
\[
2\ell+2N>\deg\sigma,
\qquad
2\ell-2+2N>\deg\tau.
\]
Then there exist
$\widehat\sigma,\widehat\tau\in\SOS_{\C}$,
both strictly positive on $\R$, such that
\[
1=w_{\mathbf e}\widehat\sigma+
\Delta_{\mathbf e}\widehat\tau,
\]
with
$ \deg\widehat\sigma=2\ell+2N $
and 
$ \deg\widehat\tau=2\ell-2+2N. $
\end{lemma}

\begin{proof}
Set
\[
(\kappa_0,\kappa_1)=
\bigl(
-\Delta_{\mathbf e}(1+x^2)^N,\,
w_{\mathbf e}(1+x^2)^N
\bigr).
\]
For $\delta>0$, define
\[
\sigma_\delta=\sigma+\delta\kappa_0,
\qquad
\tau_\delta=\tau+\delta\kappa_1.
\]
Since $(\kappa_0,\kappa_1)\in\ker\mathcal L_{\mathbf e,N}$,
we have
\[
1=w_{\mathbf e}\sigma_\delta+
\Delta_{\mathbf e}\tau_\delta.
\]
Moreover, by the choice of $N$,
$\deg\sigma_\delta=2\ell+2N$ and
$\deg\tau_\delta=2\ell-2+2N$.
We show that both $\sigma_\delta$ and $\tau_\delta$
are strictly positive on $\R$ for sufficiently small
$\delta>0$.

\medskip

By Lemma~\ref{lem:strict positive},
there exists $\varepsilon>0$ such that
$\sigma,\tau\ge\varepsilon$ on $\R$.
Since $w_{\mathbf e}>0$ on $\R$,
\[
\tau_\delta
=\tau+\delta w_{\mathbf e}(1+x^2)^N
\ge\varepsilon.
\]
On the other hand,
\[
\sigma_\delta
=\sigma-\delta\Delta_{\mathbf e}(1+x^2)^N.
\]
If $x\notin K_{\mathbf e}$, then
$\Delta_{\mathbf e}(x)<0$, so
$\sigma_\delta(x)>0$.
Since $K_{\mathbf e}$ is compact, the number
\[
\omega_N=\max_{x\in K_{\mathbf e}}
\Delta_{\mathbf e}(x)(1+x^2)^N
\]
is finite and strictly positive. Hence
\[
\sigma_\delta(x)\ge
\varepsilon-\delta\omega_N
\qquad(x\in K_{\mathbf e}).
\]
Choosing $0<\delta<\varepsilon/\omega_N$
ensures that $\sigma_\delta>0$ on $\R$.
Since every nonnegative real polynomial in one variable
belongs to $\SOS_{\C}$, we conclude that
$\sigma_\delta,\tau_\delta\in\SOS_{\C}$.
Taking $\widehat\sigma=\sigma_\delta$ and
$\widehat\tau=\tau_\delta$ completes the proof.
\end{proof}

Recall that a strictly positive real polynomial of
exact even degree $2\vd$ is an interior point of
$\SOS_{\C}\cap\R[x]_{\le2\vd}$ in $\R[x]_{\le2\vd}$.
Indeed, its quotient by $(1+x^2)^{\vd}$ is bounded
below by a positive constant on $\R$, so sufficiently
small coefficient perturbations preserve positivity.

For $N\in\mathbb N$, define
\[
\mathcal C_N=
\{(\sigma,\tau)\in V_N:\sigma,\tau\in\SOS_{\C}\}.
\]
By the preceding observation, the pair
$(\widehat\sigma,\widehat\tau)$ constructed in
Lemma~\ref{lem:strict positive plus degrees}
is an interior point of $\mathcal C_N$ in $V_N$.

\begin{lemma}\label{lem:local-stability}
Fix $\mathbf e_0\in\mathcal P_\ell$
and $N\in\mathbb N$.
Suppose there exists
$(\sigma_0,\tau_0)$ in the interior
of $\mathcal C_N$ such that
$ \mathcal L_{\mathbf e_0,N}
(\sigma_0,\tau_0)=1. $
Then there is a neighborhood $\mathcal U$ of
$\mathbf e_0$ in $\mathcal P_\ell$
such that, for every $\mathbf e\in \mathcal U$,
there exist
$(\sigma_{\mathbf e},\tau_{\mathbf e})
\in\mathcal C_N$ satisfying
$ \mathcal L_{\mathbf e,N}
(\sigma_{\mathbf e},\tau_{\mathbf e})=1. $
\end{lemma}

\begin{proof}
By Lemma~\ref{lem:L-surjective},
the map $\mathcal L_{\mathbf e_0,N}$
is surjective. Hence it admits a
right inverse
\[
\mathcal R_0:W_N\longrightarrow V_N.
\]

By Lemma~\ref{lem:h-continuity}, the coefficients of
$w_{\mathbf e}$ depend continuously on $\mathbf e$.
The same holds for $\Delta_{\mathbf e}$, since its
coefficients depend polynomially on the endpoints
by \eqref{eq:P-Q-r-Delta}.
Since $V_N$ and $W_N$ are fixed
finite-dimensional spaces,
$\mathcal L_{\mathbf e,N}$ also
depends continuously on $\mathbf e$.
Now
\[
\mathcal L_{\mathbf e_0,N}\mathcal R_0
=\operatorname{id}_{W_N}.
\]
By continuity,
$\mathcal L_{\mathbf e,N}\mathcal R_0$
remains invertible for $\mathbf e$
sufficiently close to $\mathbf e_0$.
For such $\mathbf e$, define
\[
\mathcal R_{\mathbf e}
=
\mathcal R_0
\bigl(
\mathcal L_{\mathbf e,N}\mathcal R_0
\bigr)^{-1}.
\]
Then $\mathcal R_{\mathbf e}$ is a
right inverse of $\mathcal L_{\mathbf e,N}$,
and it depends continuously on
$\mathbf e$.
Set
\[
(\sigma_{\mathbf e},\tau_{\mathbf e})
=
(\sigma_0,\tau_0)
+
\mathcal R_{\mathbf e}
\bigl(
1-\mathcal L_{\mathbf e,N}
(\sigma_0,\tau_0)
\bigr).
\]
By construction,
$ \mathcal L_{\mathbf e,N}
(\sigma_{\mathbf e},\tau_{\mathbf e})=1. $
Furthermore, since
$\mathcal L_{\mathbf e_0,N}(\sigma_0,\tau_0)=1$
and the maps $\mathcal R_{\mathbf e}$ are continuous,
\[
(\sigma_{\mathbf e},\tau_{\mathbf e})
-(\sigma_0,\tau_0)
=\mathcal R_{\mathbf e}
\bigl(1-\mathcal L_{\mathbf e,N}(\sigma_0,\tau_0)\bigr)
\longrightarrow 0
\]
as $\mathbf e\to\mathbf e_0$.
Since $(\sigma_0,\tau_0)$ is an interior
point of $\mathcal C_N$, the perturbed
pair also belongs to $\mathcal C_N$
for $\mathbf e$ sufficiently close
to $\mathbf e_0$.
Shrinking the neighborhood if necessary
proves the lemma.
\end{proof}

\begin{lemma}\label{lem:uniform-scalar-degrees}
For each fixed $\ell\ge1$ with $\mathcal P_\ell\ne\varnothing$,
there exists $N_{1,\ell}\in\mathbb N$ such that,
for every $\mathbf e\in\mathcal P_\ell$,
the identity \eqref{eq:uniform-scalar}
admits SOS multipliers satisfying
\[
\deg\sigma\le2\ell+2N_{1,\ell},
\qquad
\deg\tau\le2\ell-2+2N_{1,\ell}.
\]
\end{lemma}

\begin{proof}
Fix $\mathbf e_0\in\mathcal P_\ell$.
By Lemmas~\ref{lem:strict positive}
and~\ref{lem:strict positive plus degrees},
we can choose $N_{\mathbf e_0}\in\mathbb N$
and multipliers
$\sigma_0,\tau_0\in\SOS_{\C}$
satisfying
\[
1=w_{\mathbf e_0}\sigma_0+
\Delta_{\mathbf e_0}\tau_0,
\]
such that $(\sigma_0,\tau_0)$ is
an interior point of
$\mathcal C_{N_{\mathbf e_0}}$.

\smallskip

By Lemma~\ref{lem:local-stability},
there exists a neighborhood
$U_{\mathbf e_0}$ of $\mathbf e_0$
on which \eqref{eq:uniform-scalar}
admits SOS multipliers with degrees
at most
$ 2\ell+2N_{\mathbf e_0} $
and
$ 2\ell-2+2N_{\mathbf e_0},$
respectively.
These neighborhoods form an open cover
of $\mathcal P_\ell$.
By compactness, finitely many of them
cover $\mathcal P_\ell$.
Taking $N_{1,\ell}$ to be the largest
of the corresponding values of
$N_{\mathbf e_0}$ gives the required
uniform degree bounds.
\end{proof}

\begin{proof}[Proof of Theorem~\ref{thmC}]
Fix $\ell$ such that $\mathcal P_\ell\ne\varnothing$,
and let $N_{1,\ell}$ be as in
Lemma~\ref{lem:uniform-scalar-degrees}.
Then, for every $\mathbf e\in\mathcal P_\ell$,
there exist $\sigma,\tau\in\SOS_{\C}$ such that
\[
1=w_{\mathbf e}\sigma+\Delta_{\mathbf e}\tau,
\qquad
\deg\sigma\le2\ell+2N_{1,\ell}.
\]
Multiplying by $\Delta_{\mathbf e}$ gives
\[
\Delta_{\mathbf e}
=\Delta_{\mathbf e}^2\tau
+w_{\mathbf e}\Delta_{\mathbf e}\sigma.
\]
Set $\beta_\ell=6\ell-2+2N_{1,\ell}$.
By Theorem~\ref{thm:intervals}, for every complex
Hilbert space $\Hh$ and every hermitian
$F\in\B(\Hh)[x]$ with $\deg F\le2$, we have
\begin{equation}\label{eq:uniform-first}
F\succeq0\text{ on }K_{\mathbf e}
\quad\Longrightarrow\quad
F\in\M_{\Hh}(\{\Delta_{\mathbf e}\})_{\le\beta_\ell}.
\end{equation}

Since only finitely many values of $\ell$ occur, we set
\[
N=\max_{\ell:\,\mathcal P_\ell\ne\varnothing}\beta_\ell.
\]
By \eqref{eq:uniform-first}, every positive semidefinite hermitian quadratic operator polynomial on $K_{\mathbf e}$ admits a certificate in $\M_{\Hh}(\{\Delta_{\mathbf e}\})$ of weighted degree at most $N$, independently of $\mathbf e$ and the Hilbert space $\mathcal H$.
\medskip

Theorem~\ref{thm:quadratic-reduction}, applied separately to each $K_{\mathbf e}$, therefore gives for every $\vd\ge1$,
\[
F\succeq0\text{ on }K_{\mathbf e}
\quad\Longrightarrow\quad
F\in\M_{\Hh}(\{\Delta_{\mathbf e}\})_{\le N+2\lceil\vd/2\rceil-2}
\subseteq\M_{\Hh}(\{\Delta_{\mathbf e}\})_{\le\vd+N-1}.
\]
Thus the required degree bound holds with $C=N-1$, which depends only on $J,\mu,\omega$. For $\vd=0$, the conclusion follows by taking the square root of the positive constant operator.
\medskip

Finally, since $\Delta_K=\prod_{g\in S_K}g$, every certificate in $\M_{\Hh}(\{\Delta_K\})$ is also a certificate in $\mathcal T_{\Hh}(S_K)$ with the same weighted degree. The real Hilbert space case follows from Lemma~\ref{lem:real-descent}, and the reverse implications follow by evaluating the certificates on $K$.
\end{proof}

\section{Further Consequences}\label{sec:further}

\subsection{Positivity of operator Laurent polynomials on unions of arcs}
\label{subsec:circle}

Throughout this subsection,
arcs are closed and may be degenerate, that is, consisting of a
single point.

\medskip

Following \cite[Section~2.5]{ZalarThesis}, a finite set
\[
\mathcal S=\{g_1,\ldots,g_s\}\subseteq\C[z,z^{-1}],
\qquad
g_j\ne0,\quad g_j^*=g_j
\quad (1\le j\le s),
\]
is called a \emph{saturated description} of a nonempty
finite union of arcs $E\subseteq\mathbb T$ if
\[
E=\{z\in\mathbb T\ :\ g_j(z)\ge0,\ 1\le j\le s\},
\]
and the following conditions hold:
\begin{enumerate}[itemsep=3pt]
\item
At every nonisolated boundary point $a$ of $E$, some
generator satisfies
\[
g_j(a)=0,
\qquad g_j'(a)\ne0.
\]
\item
At every isolated point $a$ of $E$, two generators satisfy
\[
g_j(a)=g_k(a)=0,
\qquad g_j'(a)\ne0,
\qquad g_k'(a)\ne0,
\]
and $g_jg_k\le0$ on a neighborhood of $a$ in $\mathbb T$.
\end{enumerate}

Here the boundary is taken relative to $\mathbb T$,
and the derivatives are the ordinary complex derivatives
of the Laurent polynomials.

\subsubsection{Proof of Corollary~\ref{cor:circle}}
The implication from the certificate to positivity is
immediate. If $E=\mathbb T$, the assertion follows from
the operator Fej\'er--Riesz theorem
\cite{Rosenblum1968}. Assume from now on that
$E\ne\mathbb T$.

\smallskip

Set $\ve_0=0$ and $\ve_j=\deg_L g_j$ for $1\le j\le s$.
Choose distinct points $\zeta_1,\zeta_2\in\mathbb T\setminus E$
at which none of the generators $g_j$ vanishes. This is possible
because $\mathbb T\setminus E$ is nonempty and open,
whereas each nonzero Laurent polynomial has only finitely
many zeros on $\mathbb T$. For $i=1,2$, 
\[
q(x)=1+x^2,
\qquad
\lambda_i(x)=\zeta_i\frac{x-\mathrm i}{x+\mathrm i},
\qquad
K_i=\lambda_i^{-1}(E)\subset \mathbb R.
\]
 The set $K_i$ is compact, and the real polynomials
\[
\widetilde g_{ij}(x)
   =q(x)^{\ve_j}g_j(\lambda_i(x)),
\qquad 1\le j\le s,
\]
form a saturated scalar description
$
S_i
   =\{\widetilde g_{i1},\ldots,\widetilde g_{is}\}
   \subseteq\R[x]
$
of $K_i$; see \cite[Proposition~2.5.2]{ZalarThesis}.
Moreover,
$
\deg\widetilde g_{ij}=2 \ve_j,
$
since its leading coefficient is $g_j(\zeta_i)\ne0$.
We also set $\widetilde g_{i0}=1$.

\smallskip

Let $c_{S_i}$ be an additive constant obtained by
Theorem~\ref{thmA} for $S_i$, and set
$
k_i=\left\lfloor\frac{c_{S_i}}{2}\right\rfloor.
$
These integers depend only on the chosen descriptions.
Suppose that $F\succeq0$ on $E$ and $\deg_LF\le \vd$.
Then
\[
\widetilde F_i(x)=q(x)^{\vd}F(\lambda_i(x))
\]
is a hermitian operator polynomial of degree at most $2 \vd$
which is positive semidefinite on $K_i$.
Theorem~\ref{thmA} gives
\[
\widetilde F_i
 =\sum_{j=0}^{s}\widetilde g_{ij}F_{ij}^*F_{ij},
\qquad F_{ij}\in\B(\Hh)[x],
\]
with
$
2 \ve_j+2\deg F_{ij}\le2 \vd+c_{S_i}
$
whenever $F_{ij}\ne0$. Consequently,
\[
\deg F_{ij}\le \vd + k_i- \ve_j.
\]

\smallskip

To return to the circle, set
\[
\chi_i(z)=\mathrm i\frac{\zeta_i+z}{\zeta_i-z},
\qquad
h_i(z)=(z-\zeta_i)^*(z-\zeta_i).
\]
Then $\chi_i=\lambda_i^{-1}$ and
\[
q(\chi_i(z))=\frac{4}{h_i(z)}.
\]
Let
\[
B_{ij}(z)=
\begin{cases}
2^{\ve_j-\vd}(z-\zeta_i)^{\vd+k_i-\ve_j}
F_{ij}(\chi_i(z)),
& \text{if }F_{ij}\ne0,\\[2pt]
0,
& \text{if }F_{ij}=0.
\end{cases}
\] 
The degree bound on $F_{ij}$
shows that $B_{ij}\in\B(\Hh)[z]$
with
\[
\deg B_{ij}\le \vd + k_i - \ve_j.
\]
Substituting $x=\chi_i(z)$ into the certificate for
$\widetilde F_i$ and clearing denominators gives the
identity
\begin{equation}
\label{eq:proof-cor-1.1}
h_i^{k_i}F=\sum_{j=0}^{s}g_jB_{ij}^*B_{ij}.
\end{equation}

\smallskip

It remains to remove the scalar multipliers.
Since $\zeta_1\ne\zeta_2$, the hermitian Laurent polynomials
$h_1^{k_1}$ and $h_2^{k_2}$ generate the unit ideal.
They are also nonnegative on $\mathbb T$.
The positive B\'ezout identity
\cite[p.\ 38]{ZalarThesis} therefore provides
hermitian scalar Laurent polynomials $\sigma_1,\sigma_2$,
strictly positive on $\mathbb T$, such that
\[
\sigma_1h_1^{k_1}+\sigma_2h_2^{k_2}=1.
\]
By scalar Fej\'er--Riesz factorization, we may write
$\sigma_i=p_i^*p_i$ with $p_i\in\C[z]$. Thus
\[
p_1^*p_1h_1^{k_1}+p_2^*p_2h_2^{k_2}=1.
\]
The polynomials $p_i$ can be fixed independently of
$F$, $\vd$, and $\Hh$. Multiplying 
\eqref{eq:proof-cor-1.1}
by $p_i^*p_i$ and adding gives
\[
F=\sum_{j=0}^{s}g_j
   \sum_{i=1}^{2}(p_iB_{ij})^*(p_iB_{ij}).
\]
Set
\[
c_{\mathcal S}
   =\max_{i=1,2}\bigl(k_i+\deg p_i\bigr).
\]
Every nonzero factor in this expression satisfies
\[
\ve_j+\deg(p_iB_{ij})\le \vd + c_{\mathcal S}.
\]
For each $j$, the Laurent polynomial
\[
\Sigma_j
   =\sum_{i=1}^{2}(p_iB_{ij})^*(p_iB_{ij})
\]
is positive semidefinite on the whole circle and, when
nonzero, has Laurent degree at most
$\vd+c_{\mathcal S} - \ve_j$. The operator Fej\'er--Riesz theorem
therefore gives $\Sigma_j=Q_j^*Q_j$ with
$Q_j\in\B(\Hh)[z]$ and
\[
\deg Q_j\le \vd+c_{\mathcal S} - \ve_j.
\]
Taking $Q_j=0$ when $\Sigma_j=0$ completes the proof.
\qed

\subsection{Proof of Corollary \ref{cor:unbounded}}
\label{subsec:unbounded}
The assertion for $K=\R$ follows from
Proposition~\ref{prop:factorization}.
Assume therefore that $K\ne\R$.
\smallskip

Set
\[
q(x)=1+x^2,
\qquad
\lambda(x)=\frac{x-\mathrm i}{x+\mathrm i},
\qquad
\chi(z)=\mathrm i\frac{1+z}{1-z},
\qquad
h(z)=(z-1)^*(z-1).
\]
Then $\chi=\lambda^{-1}$ and
$
h(\lambda(x))=\frac{4}{q(x)}.
$
The set
\[
E=\overline{\lambda(K)}\subseteq\mathbb T
\]
is a finite union of arcs. Since every natural generator
$g_j$ has degree one or two, the expressions
\[
b_j(z)=h(z)g_j(\chi(z)),
\qquad 1\le j\le s,
\]
belong to $\C[z,z^{-1}]$ and have Laurent degree one.
Moreover,
$
\mathcal S=\{b_1,\ldots,b_s\}
$
is a saturated description of $E$; see
\cite[Proposition~2.6.3]{ZalarThesis}.
Let $c_{\mathcal S}$ be the constant supplied by
Corollary~\ref{cor:circle}, and set
$c_S=c_{\mathcal S}.$

\smallskip

Suppose that $F\succeq0$ on $K$ and $\deg F\le \vd$.
Let $m=\left\lceil\frac \vd 2\right\rceil$ and
\[
A(z)=h(z)^mF(\chi(z)).
\]
Then $A\in\B(\Hh)[z,z^{-1}]$ is hermitian,
$\deg_L A\le m$, and $A\succeq0$ on $E$.
Positivity at the point $1\in E$ follows by continuity.
Set $N=m+c_S$. Corollary~\ref{cor:circle} gives
\[
A=Q_0^*Q_0+\sum_{j=1}^{s}b_jQ_j^*Q_j,
\]
where $Q_j\in\B(\Hh)[z]$ and
$\deg Q_0\le N$,
$\deg Q_j\le N-1$ for $1\le j\le s$.
Define
\[
F_0(x)=2^{-m}(x+\mathrm i)^NQ_0(\lambda(x))
\]
and, whenever $Q_j\ne0$,
\[
F_j(x)=2^{1-m}(x+\mathrm i)^{N-1}Q_j(\lambda(x)),
\qquad 1\le j\le s,
\]
taking $F_j=0$ if $Q_j=0$.
The degree bounds on the $Q_j$ imply that
$F_j\in\B(\Hh)[x]$. Substituting $z=\lambda(x)$
into the certificate for $A$ gives
\[
q(x)^{c_S}F(x)
 =F_0(x)^*F_0(x)
  +\sum_{j=1}^{s}g_j(x) F_j(x)^*F_j(x).
\]

To obtain the asserted degree bound, choose
$\varepsilon\in\{-1,1\}$ such that $\varepsilon t\in K$
for all sufficiently large $t$.
Every nonzero summand
$
g_j(\varepsilon t)
F_j(\varepsilon t)^*F_j(\varepsilon t)
$
has a nonzero positive semidefinite leading coefficient.
Hence the terms of highest degree cannot cancel.
Consequently,
\[
\deg(g_j F_j^*F_j)
 \le\deg(q^{c_S}F)
 \le \vd+2c_S.
\]

Conversely, the certificate implies $F\succeq0$ on $K$,
since every $g_j$ is nonnegative on $K$ and $q>0$ on $\R$.
\qed

\subsection{Proof of Corollary \ref{cor:degree-preservation}}
\label{subsec:degree-preservation}
If $K$ consists of at most three points, the assertion
follows from the natural-preordering bound in
Proposition~\ref{prop:finite-set}.

\smallskip

Suppose that $K$ is any other set in the list of Corollary \ref{cor:degree-preservation}.
Then $K$ contains a nondegenerate interval.
Fix a nonzero hermitian $F\in\B(\Hh)[x]$ with
$F\succeq0$ on $K$, and write
\[
F(x)=\sum_{k=0}^{\vd}A_kx^k,
\qquad \vd=\deg F.
\]
Put $S_K=\{g_1,\ldots,g_s\}$ and
$w_\alpha=g^\alpha$ for $\alpha\in\{0,1\}^s$.
Direct the finite-dimensional subspaces $V\subseteq\Hh$
by inclusion, and let $P_V$ denote the orthogonal
projection onto $V$.
The matrix degree preservation results
\cite[Theorem 2.4.2, Corollary 2.4.4]{ZalarThesis}
and \cite[Theorems 1.3--1.5]{SunZalar}
give
\begin{equation}\label{certificate}
P_VF|_V\in\mathcal T_V(S_K)_{\le \vd}.
\end{equation}

Note that the real symmetric statements in
\cite[Theorems~1.3--1.5]{SunZalar} also give the complex
hermitian statements. Indeed, for an $n\times n$ hermitian
polynomial $H=A+iB$, with $A,B$ real, the polynomial
$
\rho(H)=\begin{pmatrix}A&-B\\ B&A\end{pmatrix}
$
is real symmetric and positive semidefinite wherever $H$ is.
Applying the real matrix result to $\rho(H)$ and then compressing its
certificate by the constant isometry
$
U=\frac1{\sqrt2}\begin{pmatrix}I_n\\-iI_n\end{pmatrix},
$
gives $U^*\rho(H)U=H.$
Each compressed sum-of-squares multiplier is a globally
positive semidefinite matrix polynomial and therefore has
a matrix Fej\'er--Riesz factorization with the same degree
bound.

\smallskip

Using \eqref{certificate} we get
\[
P_VFP_V\in\mathcal T_{\Hh}(S_K)_{\le \vd}.
\]
Moreover, for every $k$,
\[
\|P_VA_kP_V\|\le\|A_k\|,
\qquad
P_VA_kP_V\longrightarrow A_k
\quad\text{strongly}.
\]
Thus the coefficients are uniformly bounded in operator
norm and converge in the weak operator topology
to those of $F$.
Since $K$ has nonempty interior,
Lemma~\ref{lem:closed-truncations} gives
\[
F\in\mathcal T_{\Hh}(S_K)_{\le \vd}.
\]

\smallskip

For necessity, suppose that $K$ does not belong to the
stated list. We first construct a nonzero real symmetric
$2\times2$ polynomial $G$ such that
\[
G\succeq0\text{ on }K,
\qquad
G\notin\mathcal T_{\C^2}(S_K)_{\le\deg G}.
\]
We distinguish three cases.

\smallskip

\emph{Case 1: $K$ is finite.}
Write $K=\{\xi_1<\cdots<\xi_r\}$, where $r\ge4$, and put
\[
a=\xi_1,\qquad b=\xi_2,\qquad c=\xi_3,\qquad d=\xi_r,
\qquad E=\{a,b,c,d\}.
\]
Thus $K\subseteq[a,b]\cup[c,d]$.
Choose
\[
1<k<\sqrt{\frac{c-a}{b-a}}
\]
and let $G$ be the polynomial from
\cite[Example~2.4.5]{ZalarThesis}, namely
\[
\begin{split}
G(x)={}&(d-c)(x-a)(x-b)
       \begin{pmatrix}1&k\\ k&k^2\end{pmatrix}
       +(a-d)(x-b)(x-c)
       \begin{pmatrix}1&k\\ k&1\end{pmatrix}
       +(b-a)(x-c)(x-d)
       \begin{pmatrix}k^2&k\\ k&1\end{pmatrix}.
\end{split}
\]
That example proves that $G\succeq0$ on $E$ and
$G\notin\mathcal T_{\C^2}(S_E)_{\le2}$.
Since $\det G$ vanishes at $a,b,c,d$ and the leading
coefficient $L$ of $G$ satisfies
\[
\det L=(k^2-1)\bigl((b-a)(d-c)(k^2-1)
                   -(c-b)(d-a)\bigr)<0,
\]
we have
\[
\det G(x)=\det L\,(x-a)(x-b)(x-c)(x-d)\ge0
\quad\text{on }[a,b]\cup[c,d].
\]
Moreover, $G_{11}$ is strictly concave, $G_{11}(a)>0$,
and $G_{11}(d)=0$, so $G_{11}>0$ on $[a,d)$.
Together with $G(d)\succeq0$, this proves that
$G\succeq0$ on $[a,b]\cup[c,d]$, and hence on $K$.

\smallskip

If $G$ had a natural-preordering certificate on $K$ of
weighted degree at most two, every scalar weight in that
certificate would be nonnegative on $E$.
By the scalar degree-preservation theorem \cite{KMS}, each
such weight has a certificate in the natural preordering
of $E$ with weighted degree at most its own degree.
Substituting these scalar certificates into the matrix
certificate gives
$G\in\mathcal T_{\C^2}(S_E)_{\le2}$, a contradiction.

\smallskip

\emph{Case 2: $K$ is unbounded.}
Since $K$ is excluded from the list, it has either a
nondegenerate bounded interval component in addition to
an unbounded interval component, or at least two isolated points in addition to
one or two unbounded intervals.
By \cite[Theorem~3.2 and its proof]{Zalar}, there exists a
real symmetric $2\times2$ polynomial $G$ such that
$G\succeq0\text{ on }K$,
$G\notin\mathcal T_{\C^2}(S_K)$.
In particular, $G$ has no certificate of weighted degree
at most $\deg G$.

\smallskip

\emph{Case 3: $K$ is compact and infinite.}
Choose $c$ in the interior of a nondegenerate interval
component of $K$, and put
$
L=\left\{\frac1{x-c}:x\in K\setminus\{c\}\right\}.
$
The set $L$ is closed and consists of two unbounded
intervals together with the images of the remaining
components of $K$.
Since $K$ is excluded from the list, $L$ contains either
a nondegenerate bounded interval or at least two isolated
points. Thus \cite[Theorem~3.2 and its proof]{Zalar}
provides a real symmetric $2\times2$ polynomial $H$ with
$H\succeq0\text{ on }L$,
$H\notin\mathcal T_{\C^2}(S_L).$
Thus $\deg H=2m$ is even and its leading coefficient
is positive semidefinite. Therefore the polynomial
\[
G(x)=(x-c)^{2m}H\left(\frac1{x-c}\right)
\] 
is positive semidefinite
on $K$; its value at $c$ is the leading coefficient of $H$.

\smallskip

We claim that
$G\notin\mathcal T_{\C^2}(S_K)_{\le2m}$.
Otherwise, write
\[
G=\sum_\alpha w_\alpha R_\alpha,
\qquad
R_\alpha\succeq0\text{ on }\R,
\qquad
\deg(w_\alpha R_\alpha)\le2m.
\]
Let
$e_\alpha=\lceil\deg w_\alpha/2\rceil$.
Since a nonzero globally positive semidefinite matrix
polynomial has even degree,
$\deg R_\alpha\le2m-2e_\alpha$.
Consequently,
\[
\begin{split}
h_\alpha(t)&=t^{2e_\alpha}w_\alpha(c+t^{-1}),\qquad
Q_\alpha(t)=t^{2m-2e_\alpha}R_\alpha(c+t^{-1})
\end{split}
\]
are polynomials, with $h_\alpha\ge0$ on $L$ and
$Q_\alpha\succeq0$ on $\R$.
Positivity of $Q_\alpha$ at zero follows by continuity;
also, $0\notin L$ because $K$ is compact.
By scalar saturation of the natural preordering
\cite{KM}, $h_\alpha\in\mathcal T_{\C}(S_L)$.
Matrix Fej\'er--Riesz factorization of $Q_\alpha$ now gives
\[
H(t)=t^{2m}G(c+t^{-1})
    =\sum_\alpha h_\alpha(t)Q_\alpha(t)
    \in\mathcal T_{\C^2}(S_L),
\]
which is a contradiction. Since $\deg G\le2m$, the claim also
excludes a certificate of weighted degree at most $\deg G$.

\smallskip

Choose a two-dimensional subspace $V\subseteq\Hh$
and regard $G\oplus0$ as an operator polynomial on
$\Hh=V\oplus V^\perp$.
If $G\oplus0$ admitted a natural-preordering certificate
of weighted degree at most $\deg G$,
compression to $V$ would give
\[
G=\sum_\alpha w_\alpha R_\alpha,
\qquad
R_\alpha\succeq0\text{ on }\R,
\qquad
\deg(w_\alpha R_\alpha)\le \deg G,
\]
where the $R_\alpha$ are $2\times2$ matrix polynomials.
By matrix Fej\'er--Riesz factorization, each $R_\alpha$
admits a hermitian-square factorization with the same
degree bound, regarding real matrix polynomials as complex
ones when necessary. Hence
$G\in\mathcal T_{\C^2}(S_K)_{\le\deg G}$, a contradiction.
\qed

\section*{Acknowledgements}
We thank Igor Klep and Jaka Cimpri\v{c} for reading
an earlier version of this manuscript. We are grateful
to Igor Klep for asking the questions that led to
Corollaries~\ref{cor:circle} and~\ref{cor:unbounded},
and to Jaka Cimpri\v{c} for suggesting
Corollary~\ref{cor:degree-preservation}.

\section*{Declaration of AI use}
During the preparation of this manuscript, the authors used large language models (LLMs) for copy editing and for rewriting some parts of the manuscript. LLM assistance also contributed to the
development of Example~II in Subsection~\ref{ssec:exampleII} which is essential to establish Theorem~\ref{thmB}. In addition, LLMs assisted in working out
the supplementary examples in the accompanying repository~\cite{Code} by
applying methods developed in this paper. The authors have independently
verified all mathematical arguments in Example~II and in the supplementary
examples, and take full responsibility
for the content and conclusions of the manuscript and the supplementary
material.

\appendix

\section{The real Hilbert space case}\label{app:real}

The degree bounds in Theorems~\ref{thmA} and~\ref{thmC} are established for complex Hilbert spaces. In this appendix, we show that the same bounds hold for real Hilbert spaces. Recall that the operator quadratic modules and preorderings over real Hilbert spaces are defined in Subsection~\ref{ssec:Notation}, where each SOS multiplier is a sum of at most two hermitian squares.
\medskip

The complexification of $\Hh$
is the complex Hilbert space
$
\Hh_{\C}=\Hh\oplus i\Hh.
$
Every bounded real-linear operator
$A\in\B(\Hh)$ extends to a bounded complex-linear
operator $A_{\C}\in\B(\Hh_{\C})$ defined by
\[
A_{\C}(u+iv)=Au+iAv,
\qquad u,v\in\Hh.
\]
For an operator polynomial
$ F(x)=\sum_k A_kx^k\in\B(\Hh)[x], $
we define its coefficientwise complexification by
\[
F_{\C}(x) \,=\, \sum_k (A_k)_{\C} \, x^k
\in\B(\Hh_{\C})[x].
\]

\smallskip

We also recall the decomposition of a complex
operator into its real and imaginary parts.
For every $A\in\B(\Hh_{\C})$ and $\nu\in\Hh$,
write uniquely
\[
A \nu = B \nu+iD\nu,\qquad B\nu,D\nu\in\Hh.
\]
This defines bounded real-linear operators
$B,D\in\B(\Hh)$ satisfying
$ A=B_{\C}+iD_{\C}.$
We define the conjugate of $A$ by
$ \overline A=B_{\C}-iD_{\C}.$

\medskip

Every polynomial $G$ with coefficients in
$\B(\Hh_{\C})$ admits a unique decomposition
\[
G=Y_{\C}+iZ_{\C},
\]
where $Y,Z\in\B(\Hh)[x]$ satisfy
$\deg Y,\deg Z\le\deg G$, and $Y_{\C},Z_{\C}$
denote their coefficientwise complexifications.
The conjugate polynomial is then given by
$
\overline G=Y_{\C}-iZ_{\C}.
$

\medskip

\begin{lemma}\label{lem:real-descent}
Let $\Hh$ be a real Hilbert space, and let
$F\in\B(\Hh)[x]$ be hermitian.
Then, for every finite $S\subset\R[x]$ and every
nonnegative integer $n$,
\[
F\in\M_{\Hh}(S)_{\le n}
\quad\Longleftrightarrow\quad
F_{\C}\in\M_{\Hh_{\C}}(S)_{\le n}.
\]
\end{lemma}

\begin{proof}
Complexification preserves adjoints, products, and polynomial
degrees. Thus, complexifying each square factor in a real
operator certificate proves the forward implication.

\smallskip

Conversely, suppose that
\[
F_{\C} = F_0^* F_0 + \sum_{j=1}^{s}g_j F_{j}^*F_{j}, \qquad F_j\in \B(\Hh_{\mathbb C})[x], \quad j=0,1,\ldots,s,
\]
where 
$S=\{g_1,\ldots,g_s\}$, 
and the certificate
has weighted degree at most $n$.

\medskip

Let $Y_{j},Z_{j}\in\B(\Hh)[x]$ satisfy
$
F_{j}=(Y_{j})_{\C}+i(Z_{j})_{\C}.
$
Since $F_{\C}$ and the scalar generators $g_j$ are
invariant under conjugation, averaging the certificate
with its conjugate gives
\begin{align*}
F_{\C}
& \,=\, \frac12 \bigl(F_{0}^{*}F_{0} + \overline F_0^{\,*}\overline F_0 \bigr) + \frac12\sum_{j=1}^s g_j
\bigl(F_j^*F_j+
\overline F_j^{\,*}\overline F_j\bigr)\\
& \,=\, (Y_0^*Y_0+Z_0^*Z_0)_{\C} + \left(
\sum_{j=1}^s g_j
(Y_j^*Y_j+Z_j^*Z_j)
\right)_{\C}.
\end{align*}
Restricting to $\Hh$ we get
\[
F=Y_0^* Y_0+Z_0^* Z_0+
\sum_{j=1}^{s}g_j
\bigl(Y_{j}^*Y_{j}+Z_{j}^*Z_{j}\bigr),
\]
which is a real operator certificate.
Since
$\deg Y_{j},\deg Z_{j}\le\deg F_{j}$,
the resulting certificate has weighted degree at most $n$.
\end{proof}

\end{document}